\documentclass[12pt, leqno]{amsart}

\usepackage{lscape}
\usepackage{graphicx}
\usepackage{amssymb,amsmath,amsthm,amscd}
\usepackage{mathrsfs}
\usepackage{enumerate}
\usepackage[usenames,dvipsnames]{color}

\usepackage[colorlinks=true, pdfstartview=FitV, linkcolor=blue,citecolor=blue,urlcolor=blue]{hyperref}

\usepackage[all]{xy}
\allowdisplaybreaks[4]
\usepackage{oldgerm}
\usepackage{marginnote}
\usepackage{xspace}

\usepackage[normalem]{ulem}  % for strikeout

\newcommand{\nc}{\newcommand}
\numberwithin{equation}{section}

\newenvironment{blue}{\relax\color{blue}}{\hspace*{.5ex}\relax}

\newcommand{\bj}{\begin{blue}}
\newcommand{\ej}{\end{blue}}

\newcounter{mycounter}
\newcounter{myc}
\newcounter{mycounters}
\newcounter{mycountert}
\newcounter{mycounterf}
\newlength{\my}
\theoremstyle{plain}
\newtheorem{lemma}{Lemma}[section]
\newtheorem{prop}[lemma]{Proposition}
\newtheorem{theorem}[lemma]{Theorem}

\newcommand{\Prop}{\begin{prop}}
\newcommand{\enprop}{\end{prop}}
\newcommand{\Lemma}{\begin{lemma}}
\newcommand{\enlemma}{\end{lemma}}
\newcommand{\Th}{\begin{theorem}}
\newcommand{\enth}{\end{theorem}}
\newtheorem{corollary}[lemma]{Corollary}
\newcommand{\Cor}{\begin{corollary}}
\newcommand{\encor}{\end{corollary}}
\newtheorem{definition}[lemma]{Definition}
\newtheorem*{conjecture}{Conjecture}
\newcommand{\Def}{\begin{definition}}
\newcommand{\edf}{\end{definition}}
\newtheorem{sublemma}[lemma]{Sublemma}
\newcommand{\Sublemma}{\begin{sublemma}}
\newcommand{\ensub}{\end{sublemma}}

\theoremstyle{definition}
\newtheorem{remark}[lemma]{Remark}
\newtheorem{example}[lemma]{Example}
\newtheorem{Convention}[lemma]{Convention}
\newcommand{\Conv}{\begin{Convention}}
\newcommand{\enconv}{\end{Convention}}
\nc{\Con}{\begin{conjecture}}
\nc{\encon}{\end{conjecture}}
\nc{\Rem}{\begin{remark}}
\nc{\enrem}{\end{remark}}
\nc{\Ex}{\begin{example}}
\nc{\enEx}{\end{example}}

\makeatletter
\renewcommand{\p@enumii}{}
\makeatother

\newcommand{\Q}{\mathbb {Q}}

\newcommand{\Z}{{\mathbb Z}}
\newcommand{\B}{{\mathbf{B}}}
\newcommand{\A}{{\mathbf A}}

\newcommand{\D}{\mathcal{D}}

\newcommand{\R}{{\rm R}}

\newcommand{\seteq}{\mathbin{:=}}

\newcommand{\g}{{\mathfrak{g}}}

\newcommand{\Uqm}[1][{\mathfrak{g}}]{{U_q^-(#1)}}

\newcommand{\Hom}{\operatorname{Hom}}
\newcommand{\End}{\operatorname{End}}

\newcommand{\isoto}[1][]{\mathop{\xrightarrow%
[{\raisebox{.3ex}[0ex][.3ex]{$\scriptstyle{#1}$}}]%
{{\raisebox{-.6ex}[0ex][-.6ex]{$\mspace{2mu}\sim\mspace{2mu}$}}}}}
\newcommand{\Ext}{\operatorname{Ext}}

\newcommand{\eq}{\begin{eqnarray}}
\newcommand{\eneq}{\end{eqnarray}}

\newcommand{\hs}{\hspace*}

\newcommand{\To}[1][{\hs{2ex}}]{\xrightarrow{\,#1\,}}

\newcommand{\eqn}{\begin{eqnarray*}}
\newcommand{\eneqn}{\end{eqnarray*}}
\newcommand{\on}{\operatorname}

\newcommand{\bna}{\be[{\rm(a)}]}

\newcommand{\QED}{\end{proof}}
\newcommand{\Proof}{\begin{proof}}

\newcommand{\soplus}{\mathop{\mbox{\normalsize$\bigoplus$}}\limits}

\newcommand{\id}{\on{id}}
\newcommand{\ba}{\begin{array}}
\newcommand{\ea}{\end{array}}

\newcommand{\monoto}{\rightarrowtail}

\newcommand{\set}[2]{\left\{#1 \mid #2 \right\}}
\newcommand{\supp}{\operatorname{supp}}

\newcommand{\Mod}{\operatorname{Mod}}

\newcommand{\Modg}{\operatorname{{Mod}_{\mathrm{gr}}}}

\newcommand{\eqsub}{\begin{subequations}\begin{eqnarray}}
\newcommand{\eneqsub}{\end{eqnarray}\end{subequations}}

\newcommand{\ol}{\overline}

\nc{\la}{\lambda}
\nc{\lam}{\lambda}
\nc{\U}[1][\g]{U_q(#1)}
\nc{\te}{\tilde{e}}
\nc{\tei}{\tilde{e}_i}
\nc{\tf}{\tilde{f}}
\nc{\tfi}{\tilde{f}_i}
\nc{\tU}{\widetilde U_q(\g)}
\nc{\tE}{\tilde{E}}
\nc{\tF}{\widetilde{\F}}
\nc{\tK}{\widetilde{K}}

\nc{\tk}{\tilde{k}}
\nc{\tkone}{\tk_{\ol{1}}}
\nc{\teone}{\tilde{e}_{\ol{1}}}
\nc{\tfone}{\tilde{f}_{\ol{1}}}

\nc{\teibar}{\tilde{e}_{\ol{i}}} \nc{\tfibar}{\tilde{f}_{\ol{i}}}
\nc{\tki}{{\tk}_{\ol {i}}}

\nc{\BZ}{{\mathbb{Z}}}
\nc{\al}{\alpha}
\nc{\qs}{{q}}
\nc{\lan}{\langle}
\nc{\ran}{\rangle}
\nc{\re}{{\mathrm{re}}}
\nc{\wt}{\operatorname{wt}}
\nc{\ch}{\operatorname{ch}}
\nc{\Um}[1][\g]{U^-_q(#1)}
\nc{\Ue}{U^+_q(\g)}
\nc{\eps}{\varepsilon}
\nc{\vphi}{\varphi}
\nc{\sphi}{\varphi^*}
\nc{\seps}{\varepsilon^*}

\nc{\nn}{\nonumber}

\nc{\vp}{\varpi}
\nc{\cls}{{\operatorname{cl}}}
\nc{\Wt}{{\operatorname{Wt}}}
\nc{\Us}{U'_q(\g)}
\nc{\La}{\Lambda}
\nc{\tLa}{\widetilde\Lambda}
\nc{\ro}{{\rm(}}
\nc{\rf}{{\rm)}}
\nc{\norm}{{\mathrm{norm}}}
\nc{\qbox}{\quad\mbox}
\nc{\braid}{{\mathfrak{B}}}
\nc{\Ad}{\operatorname{Ad}}
\nc{\Aut}{\operatorname{Aut}}
\nc{\dt}[1]{\tilde{\tilde #1}}
\nc{\Sn}{S^{{\mathrm{norm}}}}
\nc{\aff}{{\rm{aff}}}
\nc{\rk}{{\mathrm{rk}}}
\nc{\tP}{\widetilde{P}}
\nc{\tW}{\widetilde{W}}
\nc{\Dyn}{\mathrm{Dyn}}
\nc{\tD}{\widetilde{\Delta}}
\nc{\height}[1]{{\operatorname{ht}}(#1)}
\nc{\bl}{\bigl(}
\nc{\br}{\bigr)}
\nc{\Hecke}{\mathrm{H}}
\nc{\HA}{\Hecke^{\mathrm{A}}}
\nc{\HB}{\Hecke^{\mathrm{B}}}
\newcommand{\scbul}{{\,\raise1pt\hbox{$\scriptscriptstyle\bullet$}\,}}
\nc{\vac}{{\phi}}
\nc{\Bt}{\B_\theta(\g)}
\nc{\be}{\begin{enumerate}}
\nc{\ee}{\end{enumerate}}
\nc{\low}{{\mathrm{low}}}
\nc{\upper}{{\mathrm{up}}}
\nc{\Zodd}{\Z_{\mathrm{odd}}}
\nc{\Ft}[1][n]{\mathbb{P}\mathrm{ol}_{#1}}
\nc{\Ftf}[1][n]{\widetilde{\mathbb{P}\mathrm{ol}}_{#1}}
\nc{\KA}{\on{K}^{\mathrm{A}}}
\nc{\KB}{\on{K}^{\mathrm{B}}}
\nc{\Res}{\on{Res}}

\nc{\tphi}{\tilde{\varphi}}
\nc{\CO}{\mathscr{O}}
\nc{\inte}{\mathrm{int}}
\nc{\Oint}{\mathcal{O}^{\ge0}_{\inte}}
\nc{\vs}{\vspace*}
\nc{\tLt}{\widetilde{L}}
\nc{\tL}{\widetilde{\Lambda}}
\nc{\tu}{\tilde{u}}
\nc{\noi}{\noindent}
\nc{\heigh}{\mathfrak{t}}
\nc{\lowest}{\mathfrak{l}}
\nc{\rootl}{\mathsf{Q}}
\nc{\cl}{\colon}
\nc{\uqpg}{U'_q(\mathfrak g)}
\nc{\uq}{\uqpg}
\nc{\Oh}{\widehat{\mathcal{O}}}

\nc{\KLR}{KLR algebra}
\nc{\KLRs}{KLR algebras}
\nc{\cor}{\mathbf{k}}
\nc{\cora}{{\cor(A)}}
\nc{\haut}{\mathrm{ht}}
\nc{\tens}{\mathop\otimes}
\nc{\gmod}{\mbox{-$\mathrm{gmod}$}}
\nc{\gMod}{\mbox{-$\mathrm{gMod}$}}
\nc{\proj}{\mbox{-$\mathrm{proj}$}}
\nc{\gproj}{\mbox{-$\mathrm{gproj}$}}
\nc{\smod}{\mbox{-$\mathrm{mod}$}}
\nc{\h}{\mathfrak h}
\nc{\Rnorm}{R^{\rm{norm}}}
\newcommand{\F}{\mathfrak{F}}   % Functor
\nc{\Vhat}{\widehat{V}}
\def\T{{\mathcal T}}

\nc{\fd}[1][A]{\on{\mathrm{flat.dim}_{#1}}}
\nc{\bP}{{\mathbb{P}}}
\nc{\bPh}{\widehat{\mathbb{P}}}
\nc{\bK}[1][{n}]{\widehat{\mathbb{K}}_{#1}}
\nc{\bV}[1][{n}]{\widehat{V}^{\otimes{#1}}}
\nc{\bVK}[1][{n}]{\widehat{V}^{\otimes{#1}}_{\widehat{\mathbb{K}}}}
\nc{\hV}{\widehat{V}}
\nc{\opp}{\mathrm{opp}}
\nc{\col}{\colon}
\nc{\bnum}{\be[{\rm(i)}]}
\nc{\oep}{\epsilon}
\nc{\qtext}{\quad\text}
\nc{\qtextq}[1]{\quad\text{#1}\quad}
\nc{\longtwoheadrightarrow}[1][]{\xymatrix{\ar@{->>}[r]^-{{#1}}&}}
\nc{\epiTo}[1][]{\longtwoheadrightarrow[{#1}]}
\nc{\epito}{\twoheadrightarrow}
\nc{\monoTo}[1][]{\xymatrix{\ar@{>->}[r]^-{{#1}}&}}
\nc{\sym}{\mathfrak{S}}
\nc{\inp}[1]{{({#1})_{\mathrm{n}}}}
\nc{\rtl}{\rootl}
\nc{\wtd}{\widetilde}
\nc{\etens}{\boxtimes}
\nc{\ds}[1]{\mathrm{d}(#1)}
\nc{\rmat}[1]{{\mathbf{r}}_%
{\mspace{-2mu}\raisebox{-.6ex}{${\scriptstyle{#1}}$}}}
\nc{\rmats}[1]{{\mathbf{r}}_%
{\mspace{-2mu}\raisebox{-.6ex}{${\scriptscriptstyle{#1}}$}}}
\nc{\shc}{\mathcal{C}}
\nc{\shs}{\mathcal{S}}
\nc{\Fct}{{\on{Fct}}}
\nc{\tC}{\widetilde{\shc}}
\nc{\Zp}{\Z_{\ge0}}
\nc{\tPhi}{\widetilde{\Phi}}
\nc{\tT}{{\widetilde{\T}}}
\nc{\Ob}{\on{Ob}}
\nc{\Img}{\on{Im}}
\nc{\Ab}{\mathcal{A}^{\mathrm{big}}}
\nc{\Sb}{\mathcal{S}^{\mathrm{big}}}
\nc{\As}{\mathcal{A}}
\nc{\Ss}{\mathcal{S}}
\nc{\ntens}{\widetilde{\otimes}}
\nc{\hR}{\widehat{R}}
\nc{\ts}{\tilde{s}}
\nc{\sho}{\mathcal{O}}
\nc{\bc}{\begin{cases}}
\nc{\ec}{\end{cases}}
\nc{\slnh}{{\widehat{\mathfrak{sl}}_N}}
\nc{\UA}{U_q'(\slnh)}
\nc{\KR}{R_K}
\nc{\cQ}{\mathcal{Q}}
\nc{\Irr}{\mathcal{I}rr}
\nc{\tQ}{\widetilde{\cQ}}
\nc{\bs}{\mathbf{s}}
\nc{\bL}{\mathbb{L}}
\nc{\tg}{\tilde{g}}

\nc{\conv}{{\mathbin{\scalebox{1.1}{$\mspace{1.5mu}\circ\mspace{1.5mu}$}}}}
\nc{\shconv}{\mathbin{\large\diamond}}
\nc{\hconv}{{\scalebox{1.2}{$\mathbin\diamond$}}}

\nc{\Rm}{R^{\mathrm{norm}}}

\renewcommand{\Im}{\on{Im}}

\nc{\de}{\on{\textfrak{d}}}

\nc{\xmono}{\ar@{>->}}
\nc{\xepi}{\ar@{->>}}
\nc{\db}[1]{\raisebox{-.5ex}[2ex][1.8ex]{$#1$}}
\nc{\wb}[1]{\mbox{$\rule[-1.1ex]{0ex}{2ex}#1$}}
\nc{\univ}{\mathrm{univ}}
\nc{\rM}{{}^*\mspace{-2mu}M}
\nc{\lM}{M^*}
\nc{\uqm}{\uq\smod}
\nc{\tR}{\widetilde{R}_{\gamma,\beta}}
\nc{\tx}{\tilde{x}}
\nc{\bi}{\mathbf{i}}
\nc{\ttau}{\widetilde{\tau}}

\nc{\tEnd}{\on{E\textsc{nd}}}
\nc{\tHom}{\on{H\textsc{om}}}

\nc{\K}{{J}}
\nc{\Kex}{{\K}_{\mathrm{ex}}}
\nc{\Kfr}{{\K}_{\mathrm{f\mspace{.01mu}r}}}
\nc{\coro}{\cor}
\nc{\tB}{\widetilde{B}}
\nc{\seed}{\mathscr{S}}
\nc{\up}{\mathrm{up}}
\nc{\bfa}{\mathbf{a}}

\nc{\bM}{\bar{M}}
\nc{\bN}{\bar{N}}
\nc{\Ma}{\mathsf{M}}
\nc{\Na}{\mathsf{N}}
\nc{\Ht}[1]{\mathrm{ht}(#1)}
\nc{\ap}[1][i]{\textfrak{p}_{#1}}
\nc{\pM}{\pi_\Ma}
\nc{\dM}[1][{\Ma}]{d_{#1}}
\nc{\dN}{d_{\Na}}
\nc{\z}[1][{\Ma}]{z_{#1}}
\nc{\zN}{z_\Na}
\nc{\rlQ}{\rtl}
\newcommand{\cmA}{\mathsf{A}}  % Cartan matrix
\newcommand{\wlP}{\mathsf{P}}   % weight lattice
\newcommand{\weyl}{\mathsf{W}}  % Weyl group
\newcommand{\pr}{\Phi_+}
\nc{\bR}{\cor}
\newcommand{\sg}{\mathfrak{S}}   % symmetric group
\nc{\Po}{\wlP}
\nc{\rtlp}{\rtl_+}
\nc{\qQ}{\mathcal{Q}}
\nc{\bQ}{\ol{\qQ}}
\nc{\qQD}{\qQ^{\mathcal{D}}}
\newcommand{\Fa}{ (\mathcal{F}\text{\rm-1})  }
\newcommand{\Fb}{ (\mathcal{F}\text{\rm-2})  }
\newcommand{\Fc}{ (\mathcal{F}\text{\rm-3})  }
\newcommand{\Fd}{ (\mathcal{F}\text{\rm-4})  }
\newcommand{\Fe}{ (\mathcal{F}\text{\rm-5})  }
\nc{\zj}[1][j]{\mathsf{z}_{#1}}
\nc{\Rmat}{\mathsf{R}}
\nc{\alD}{\al^{\mathcal{D}}}
\nc{\hD}{h^{\mathcal{D}}}
\nc{\rtlD}{\rtl^{\mathcal{D}}}
\nc{\rtlpD}{\rtl_+^{\mathcal{D}}}
\nc{\PoD}{\Po^{\mathcal{D}}}
\nc{\FD}{\F^{\mathcal{D}}}
\nc{\fmod}{\mbox{-$\mathrm{fmod}$}}
\nc{\Modf}{\on{Mod_{fg}}}
\nc{\rr}[1][j]{\mathsf{r}_{#1}}
\nc{\KM}{\mathrm{K}}
\nc{\zr}{0}
\nc{\Cv}{\mathrm{C}}
\nc{\RD}{R^\mathcal{D}}
\nc{\fc}{\mathbf{c}}
\nc\Aq[1][{\g^+}]{A_q(#1)}
\nc{\Uam}[1][\g]{U_{\A}^-(#1)}
\nc{\xd}{x^\mathcal{D}}
\nc{\td}{\tau^\mathcal{D}}
\nc{\LD}{L^\D}
\nc{\prolim}[1][]{\mathop{\varprojlim}\limits_{#1}}

\begin{document}

\title[Affinizations and R-matrices for quiver Hecke algebras] {Affinizations and R-matrices for quiver Hecke algebras}

\author[Masaki Kashiwara]{Masaki Kashiwara}
\thanks{The research of the first author was supported by Grant-in-Aid for Scientific Research (B) 22340005, Japan Society for the Promotion of Science.}
\address[Masaki Kashiwara]{Research Institute for Mathematical Sciences, Kyoto University,
Kyoto 606-8502, Japan \& Korea Institute for Advanced Study, Seoul 02455, Korea  }
\email[Masaki Kashiwara]{masaki@kurims.kyoto-u.ac.jp}

\author[Euiyong Park]{Euiyong Park}
\thanks{The research of the second author was supported by the National Research Foundation of Korea(NRF) Grant funded by the Korean Government(MSIP)(NRF-2014R1A1A1002178).}
\address[Euiyong Park]{Department of Mathematics, University of Seoul, Seoul 02504, Korea}
\email[Euiyong Park]{epark@uos.ac.kr}


\date{February 9, 2018}

\begin{abstract}
We introduce the notion of affinizations and R-matrices for arbitrary quiver Hecke algebras.
It is shown that they enjoy similar properties to those
for symmetric quiver Hecke algebras.
We next define a duality datum $\mathcal{D}$ and construct
a tensor functor $\F^{\mathcal{D}}\cl \Modg(R^{\mathcal{D}})\to\Modg(R) $ between graded module categories of quiver Hecke algebras $R$ and $R^{\mathcal{D}}$
 arising from $\mathcal{D}$.
The functor $\F^{\mathcal{D}}$ sends finite-dimensional modules to finite-dimensional modules, and is exact when $R^{\mathcal{D}}$ is of finite type.
It is proved that affinizations of real simple modules and their R-matrices give a duality datum. Moreover, the corresponding duality functor sends a simple module to a simple module or zero when $R^{\mathcal{D}}$ is of finite type.
We give several examples of the functors $\F^{\mathcal{D}}$ from the graded module category of the quiver Hecke algebra of type
$D_\ell$, $C_\ell$, $B_{\ell-1}$, $A_{\ell-1}$
to that of type $A_\ell$, $A_\ell$, $B_\ell$, $B_\ell$, respectively.
\end{abstract}

\maketitle

%\tableofcontents

\vskip 2em

\section*{Introduction}

The {\it quiver Hecke algebras} (or {\it Khovanov-Lauda-Rouquier algebras}), introduced by Khovanov-Lauda (\cite{KL09, KL11}) and Rouquier (\cite{R08}) independently, are $\Z$-graded algebras which provide a categorification for the negative half of a quantum group.
The algebras are a vast generalization of affine Hecke algebras of type $A$ in the direction of categorification (\cite{BK09, R08}), and they have special graded quotients, called {\it cyclotomic quiver Hecke algebras},
 which categorify irreducible integrable highest weight modules (\cite{KK11}).
When the quiver Hecke algebras are {\it symmetric},
we can study them more deeply.

%\vs{1ex}
\noi\addtolength{\my}{-2ex}
$\bullet$\
\parbox[t]{\my}%
{First of all, it is known that
the upper global basis corresponds to the set of the isomorphism classes of simple modules over symmetric quiver Hecke algebras (\cite{R11, VV09}).}

\vs{.5ex}
\noi
$\bullet$\
\parbox[t]{\my}{The KLR-type quantum affine Schur-Weyl duality functor
was constructed in \cite{K^3} using symmetric quiver Hecke algebras
and R-matrices of quantum affine algebras.
This functor has been studied in various types (\cite{K^3b,KKKO14c,KKKO15b}).}
\setlength{\my}{\textwidth}

The notion of {\it R-matrices} for symmetric quiver Hecke algebras was introduced in \cite{K^3}. The R-matrices are special homomorphisms defined by using intertwiners and affinizations.
It turned out that the R-matrices have very good properties with {\it real} simple modules (\cite{KKKO14}).
They also had taken an important role as a main tool in studying a monoidal categorification of quantum cluster algebras (\cite{KKKO18}).

Let us explain the construction of R-matrices in \cite{K^3} briefly.
We assume that the quiver Hecke algebra $R$ is symmetric. Let $M$ be an $R$-module and $M_z$ its affinization.
The $R$-module $M_z$ is isomorphic to $\bR[z]\otimes_\bR M$ as a $\bR$-vector space. The actions of $e(\nu)$ and $\tau_i$ on $M_z$ are the same as those on $M$, but
the action of $x_i$ on $M_z$ is equal to
the action $x_i$ on $M$ added by $z$ (see $\eqref{Eq: affinization1}$).
For $R$-modules $M$ and $N$, we next consider the homomorphism $R_{M_z, N_{z'}} \in \tHom_R(M_z \conv N_{z'},  N_{z'} \conv M_z ) $ given by using intertwiners (see $\eqref{Eq: RMN}$).
Here $\tHom$ denotes the non-graded homomorphism space (see \eqref{eq:tHom}).
We set
$$\Rm_{M_z, N_{z'}} := (z'-z)^{-s} R_{M_z, N_{z'}}, \quad \mathbf{r}_{M,N} := \Rm_{M_z, N_{z'}}|_{z=z'=0},$$
where $s$ is the order of the zero of $R_{M_z, N_{z'}}$. Then the morphisms $\Rm_{M_z, N_{z'}}$ and $\mathbf{r}_{M,N}$ are non-zero, commute with the spectral parameters $z$, $z'$, and satisfy the braid relations.
Here,  in defining $M_z$ and $\mathbf{r}_{M,N}$,
we crucially use the fact that $R$ is symmetric.

In this paper, we introduce and investigate the notion of
affinizations and R-matrices for {\em arbitrary} quiver Hecke algebras, and construct a new duality functor between finitely generated graded module categories of quiver Hecke algebras.
The affinizations defined in this paper generalize the affinizations $M_z$ for symmetric quiver Hecke algebras. The root modules given in \cite{BKM14} are examples of affinizations.

We then define a tensor functor $\F^{\mathcal{D}}: \Modg(R^{\mathcal{D}})\To\Modg(R) $ between graded module categories of quiver Hecke algebras $R$ and $R^{\mathcal{D}}$,
which arises from a {\it duality datum} $\mathcal{D}$ consisting of certain $R$-modules and their homomorphisms. This is inspired by the KLR-type quantum affine Schur-Weyl duality functor in \cite{K^3}.
The functor $\F^{\mathcal{D}}$ sends finite-dimensional modules to finite-dimensional modules. It is exact when $R^{\mathcal{D}}$ is of finite type.
We show that affinizations of real simple modules and their R-matrices give a duality datum.
The corresponding duality functor sends a simple module to a simple module or zero when $R^{\mathcal{D}}$ is of finite type.

Here is a brief description of our work. Let $R(\beta)$ be an arbitrary quiver Hecke algebra. We define an affinization $(\Ma, \z)$ of a simple $R(\beta)$-module $\bM$ to be an $R(\beta)$-module $\Ma$ with a homogeneous endomorphism $\z \in \End_R(\Ma)$ and an isomorphism $\Ma/\z\Ma \simeq \bM$ satisfying the conditions in Definition \ref{def:aff}.

We then study the endomorphism rings of affinizations and the
homomorphism spaces between convolution products of  simple modules and their affinizations.
For a non-zero $R$-module $N$,
let $s$ be the largest integer such that $R_{\Ma, N} (\Ma \conv N ) \subset \z^s N \conv \Ma$.
We set
$$
\Rm_{\Ma,N} =\z^{-s}R_{\Ma,N}\cl \Ma\conv N\to N\conv \Ma,
$$
and denote by $ \rmat{\bM,N}\cl \bM\conv N\to N\conv \bM $ the homomorphism induced by $\Rm_{\Ma,N}$.
By the definition $\rmat{\bM,N}$ never vanishes. The R-matrix
$\rmat{\bM,N}$ has similar properties to R-matrices for symmetric quiver Hecke algebras (Proposition \ref{Prop: simple-unique}).
Proposition \ref{prop:MNNM} tells us that, if
$(\Ma,\z)$ and $(\Na,\z[\Na])$ are affinization of simple modules $\bM$ and $\bN$ and one of $ \bM$ and $\bN$ is real (see \eqref{hyp:MN}), then
\bnum
\item $\tHom_{R[\z,\zN]}(\Ma\conv \Na,\Ma\conv \Na) =\cor[\z,\zN]\id_{\Ma\circ \Na},$
\item $\tHom_{R[\z,\zN]}(\Ma\conv \Na,\Na\conv \Ma)$ is a free $\cor[\z,\zN]$-module of rank one.
\ee
Here, $\tHom$ denotes the space of non-graded homomorphisms (see \eqref{eq:tHom}).
We define $\Rm_{\Ma,\Na}$ as a generator of the $\cor[\z,\zN]$-module $\tHom_{R[\z,\zN]}(\Ma\conv \Na,\Na\conv \Ma)$.
Then $\Rm_{\Ma,\Na}$ commutes with $\z$ and $\zN$ by the construction, and
we prove that
$\Rm_{\Ma,\Na}|_{\z=\zN=0} \in\tHom(\bM\conv\bN,\bN\conv\bM)$
does not vanish and coincides with $\rmat{\bM, \bN}$
up to a constant multiple  in Theorem \ref{Thm: specializetion of Rm}.

We next define the duality datum $\mathcal{D}=\{\beta_j, M_j, \zj,\rr,\Rmat_{j,k}\}_{j,k\in J}$ axiomatically. Here,
$J$ is a finite index set, and
\begin{align*}
&M_j \in \Modg(R(\beta_j)), \quad \qquad \mathsf{z}_j \in \tEnd_{R(\beta_j)}(M_j),  \\
&\mathsf{r}_j \in \tEnd_{R(2\beta_j)}(M_j \circ  M_j),\quad \mathsf{R}_{j,k} \in \tHom_{R(\beta_j + \beta_k )}(M_j \circ  M_k, M_k \circ  M_j),
\end{align*}
satisfying certain conditions given in Definition \ref{Def: duality}.
We construct the generalized Cartan matrix $\cmA^{\mathcal{D}}$ and the polynomial parameters $\qQD_{i,j}(u,v)$ out of the duality datum $\mathcal{D}$
and consider the quiver Hecke algebra $R^{\mathcal{D}}$ corresponding to $\cmA^{\mathcal{D}}$ and $\qQD_{i,j}(u,v)$.
For $\gamma \in {\rlQ}_+^{\mathcal{D}}$ with $m = \Ht{\gamma} $, we define
$$ \Delta^{\mathcal{D}}( \gamma )\seteq\soplus_{\mu \in J^{ \gamma}} \Delta^{\mathcal{D}}_{\;\mu}, $$
where
$$ \Delta^{\mathcal{D}}_{\;\mu} \seteq M_{\mu_1} \conv M_{\mu_2} \conv \cdots \conv
M_{\mu_m} \ \text{for $\mu = (\mu_1, \mu_2, \ldots, \mu_m) \in J^{\gamma}$.} $$
It turns out that $\Delta^{\mathcal{D}}( \gamma )$ has an $(R,  R^{\mathcal{D}})$-bimodule structure (Theorem \ref{Thm: bimodule}), and
 we obtain the duality functor $\F^{\mathcal{D}}\cl \Modg(R^{\mathcal{D}})\to\Modg(R) $
by tensoring $\Delta^{\mathcal{D}}( \gamma )$. Theorem \ref{Thm: main thm} tells that $\F^{\mathcal{D}}$ is a tensor functor and sends finite-dimensional modules to finite-dimensional modules.
Moreover it is exact when $\cmA^{\mathcal{D}}$ is of finite type.
Affinizations of real simple modules and their R-matrices provide
a duality functor which enjoys extra good properties (Theorem \ref{Thm: duality and affinization}).

Several examples of duality functors $\F^{\mathcal{D}}$ are given in Sections
\ref{Sec: Ex} and \ref{Sec: Ex2}.
In Example \ref{Ex: D}, we construct a duality functor $\F^{\mathcal{D}}$ from the graded module category of the quiver Hecke algebra of type $D_\ell$ to that of type $A_\ell$.
The other examples are in non-symmetric cases.
We discuss a duality functor from type $C_\ell$ to type $A_\ell$
in Example \ref{Ex: C},
 and ones from types $B_{\ell-1}$ and $A_{\ell-1}$ to type $B_\ell$
in Examples \ref{Ex: B1} and \ref{Ex: B2}.

\vs{1em}
\noindent
{\bf Acknowledgments.}
The authors would like to thank Myungho Kim for fruitful discussions.

\vskip 2em

\section{Preliminaries}

\subsection{Quantum groups} Let $I$ be an index set.
\Def\label{def:cartan}
A Cartan datum is a  quintuple  $\bl \cmA,\wlP,\Pi,\Pi^\vee,(\cdot,\cdot)\br$
consisting of
\bna
\item a free abelian group $\wlP$, called the {\em weight lattice},
\item $\Pi = \set{\alpha_i}{i\in I} \subset \wlP$,
called the set of {\em simple roots},
\item $\Pi^{\vee} = \set{ h_i}{i\in I} \subset \wlP^{\vee}\seteq
\Hom( \wlP, \Z )$, called the set of {\em simple coroots},
\item \label{item:inner}
a $\Q$-valued symmetric bilinear form $(\cdot\, ,\, \cdot)$  on $\Po$,
which satisfy
\be[{\rm(1)}]
\item  $(\al_i,\al_i)\in 2\Z_{>0}$ for any $i\in I$, \label{hyp:even}
\item $\lan h_i,\la\ran=\dfrac{2(\al_i,\la)}{(\al_i,\al_i)}$ for any $i\in I$and $\la\in\Po$, \label{item:h}
\item $\cmA \seteq(\lan h_i,\al_j\ran)_{i,j\in I}$ is
a {\em generalized Cartan matrix}, i.e.,
$\lan h_i,\al_i\ran=2$ for any $i\in I$ and
$\lan h_i,\al_j\ran\in\Z_{\le0}$ if $i\not=j$,
\item $\Pi$ is a linearly independent set,
\item for each $i\in I$, there exists $\Lambda_i \in \wlP$
such that $\langle h_j, \Lambda_i \rangle = \delta_{ij}$ for any $j\in I$.
\ee
\ee
\edf
Let us write $\rlQ = \bigoplus_{i \in I} \Z \alpha_i$ and $\rlQ_+ = \sum_{i\in I} \Z_{\ge 0} \alpha_i$.
For $\beta=\sum_{i \in I} k_i \alpha_i \in \rlQ_+$, set $\Ht{\beta}
=\sum_{i \in I} k_i$.
The {\em Weyl group} $\weyl$ associated with the Cartan datum
is the subgroup of $\Aut(\wlP)$ generated by
the reflections $\{r_i\}_{i\in I}$ defined by
$$r_i(\lambda)=\lambda-\langle h_i, \lambda\rangle\alpha_i \qquad
\text{for $\lambda\in \wlP$.}$$

Let $\g$ be the Kac-Moody algebra associated with a Cartan datum
$(\cmA,\wlP, \Pi, \Pi^\vee, (\cdot,\cdot))$
and $\pr$ the set of positive roots of $\g$.
We denote by $U_q(\g)$ the corresponding quantum group, which is an associative algebra over $\Q(q)$ generated by $e_i$, $f_i$ $(i\in I)$ and $q^h$
$(h\in \wlP^\vee)$ with certain defining relations (see \cite[Chap.\ 3]{HK02} for details).
Set $\A=\Z[q,q^{-1}]$. We denote by $U_\A^-(\g)$ the subalgebra of $U_q(\g)$ generated by $f_i^{(n)} := f_i^n / [n]_i!$ for $i\in I$ and $n\in \Z_{\ge0}$,
where  $q_i = q^{ (\alpha_i, \alpha_i)/2 }$ and
\begin{equation*}
 \begin{aligned}
 \ &[n]_i =\frac{ q^n_{i} - q^{-n}_{i} }{ q_{i} - q^{-1}_{i} },
 \ &[n]_i! = \prod^{n}_{k=1} [k]_i.
 \end{aligned}
\end{equation*}

\subsection{Quiver Hecke algebras}\
Let $\bR$ be a field. For $i,j\in I$, we take polynomials
$\qQ_{i,j}(u,v) \in \bR[u,v]$ such that

\noi
(i) $\qQ_{i,j}(u,v) = \qQ_{j,i}(v,u)$, and

\noi
(ii) it is of the form
\begin{align} \label{Eq: Qij}
\qQ_{i,j}(u,v) = \left\{
                 \begin{array}{ll}
                   \sum_{2(\alpha_i , \alpha_j) + p(\alpha_i , \alpha_i) + q(\alpha_j , \alpha_j) = 0} t_{i,j;p,q} u^pv^q & \hbox{if } i \ne j,\\
                   0 & \hbox{if } i=j,
                 \end{array}
               \right.
\end{align}
where $t_{i,j;-a_{ij},0} \in  \bR^{\times}$.
We set
\eq
&&\bQ_{i,j}(u,v,w)=\dfrac{\qQ_{i,j}(u,v)-\qQ_{i,j}(w,v)}{u-w}\in\cor[u,v,w].
\eneq
For $\beta\in\rtl_+$ with $\Ht{\beta}=n$, set
$$I^\beta\seteq\set{\nu=(\nu_1, \ldots, \nu_n ) \in I^n}
{\sum_{k=1}^n\alpha_{\nu_k} = \beta}.$$
The symmetric group $\mathfrak{S}_n = \langle s_k \mid k=1, \ldots, n-1 \rangle$ acts on $I^\beta$ by place permutations.

\Def \ For $\beta\in\rtlp$,
the {\em quiver Hecke algebra} $R(\beta)$ associated with $\cmA$ and $(\qQ_{i,j}(u,v))_{i,j\in I}$
is the $\bR$-algebra generated by
$$\set{e(\nu)}{\nu \in I^\beta}, \;\set{x_k}{1 \le k \le n},
 \;\set{\tau_l}{1 \le l \le n-1} $$
satisfying the following defining relations:

\begin{equation} \label{Eq: defining relations}
\begin{aligned}
& e(\nu) e(\nu') = \delta_{\nu,\nu'} e(\nu),\ \sum_{\nu \in I^{\beta}} e(\nu)=1,\
x_k e(\nu) =  e(\nu) x_k, \  x_k x_l = x_l x_k,\\
& \tau_l e(\nu) = e(s_l(\nu)) \tau_l,\  \tau_k \tau_l = \tau_l \tau_k \text{ if } |k - l| > 1, \\[5pt]
&  \tau_k^2 e(\nu) = \qQ_{\nu_k, \nu_{k+1}}(x_k, x_{k+1}) e(\nu), \\[5pt]
&  (\tau_k x_l - x_{s_k(l)} \tau_k ) e(\nu) = \left\{
                                                           \begin{array}{ll}
                                                             -  e(\nu) & \hbox{if } l=k \text{ and } \nu_k = \nu_{k+1}, \\
                                                               e(\nu) & \hbox{if } l = k+1 \text{ and } \nu_k = \nu_{k+1},  \\
                                                             0 & \hbox{otherwise,}
                                                           \end{array}
                                                         \right. \\[5pt]
&( \tau_{k+1} \tau_{k} \tau_{k+1} - \tau_{k} \tau_{k+1} \tau_{k} )  e(\nu) \\[4pt]
&\qquad \qquad \qquad = \left\{
                                                                                   \begin{array}{ll}
\bQ_{\,\nu_k,\nu_{k+1}}(x_k,x_{k+1},x_{k+2}) e(\nu) & \hbox{if } \nu_k = \nu_{k+2}, \\
0 & \hbox{otherwise}. \end{array}
\right.\\[5pt]
\end{aligned}
\end{equation}
\edf

The algebra $R(\beta)$ has the $\Z$-graded algebra structure given by
\begin{align} \label{Eq: grading}
\deg(e(\nu))=0, \quad \deg(x_k e(\nu))= ( \alpha_{\nu_k} ,\alpha_{\nu_k}), \quad  \deg(\tau_l e(\nu))= -(\alpha_{\nu_{l}} , \alpha_{\nu_{l+1}}).
\end{align}

For $\beta\in \rtl_+$, let us denote by $\Mod(R(\beta))$ the category of $R(\beta)$-modules and by $R(\beta)\smod$ the category of finite-dimensional $R(\beta)$-modules.

We denote by $\Modg\bl R(\beta)\br$ the category of graded $R(\beta)$-modules
and by $R(\beta)\gmod$ the category of finite-dimensional graded
$R(\beta)$-modules.
We denote by $\Modf\bl R(\beta)\br$ the full subcategory of
$\Modg\bl R(\beta)\br$ consisting of finitely generated graded $R(\beta)$-modules.
Their morphisms are homogeneous of degree zero.
Hence, $\Mod\bl R(\beta)\br$, $R(\beta)\smod$, $\Modg\bl R(\beta)\br$,
$R(\beta)\gmod$ and $\Modf\bl R(\beta)\br$ are abelian categories.
We set $\Modg(R)\seteq\soplus_{\beta\in\rtl^+}\Modg\bl R(\beta)\br$,
$R\smod\seteq\soplus_{\beta\in\rtl^+} R(\beta)\smod$, etc.
The objects of $\Modg(R)$ are sometimes simply called $R$-modules.

We denote by $R(\beta)\text{-proj}$
the full subcategory of $\Modg(R(\beta))$
consisting of finitely generated projective graded $R(\beta)$-modules.

Let us denote by $q$ the {\em grading shift functor},
i.e., $(qM)_k = M_{k-1}$ for a graded module $M = \bigoplus_{k \in \Z} M_k $.

For $\nu\in I^\beta$ and $\nu'\in I^{\beta'}$, let $e(\nu, \nu')$ be the idempotent corresponding to the concatenation
$\nu\ast\nu'$ of
$\nu$ and $\nu'$, and set
$$
e(\beta, \beta') := \sum_{\nu \in I^\beta, \nu' \in I^{\beta'}} e(\nu, \nu').
$$
For an $R(\beta)$-module $M$ and an $R(\beta')$-module $N$, we define an $R(\beta+\beta')$-module $M \conv N$ by
$$ M \conv N := R(\beta+\beta') e(\beta, \beta') \otimes_{R(\beta) \otimes R(\beta')} (M \otimes N). $$

We denote by $M \hconv N$ the head of $M \conv N$.

For a graded $R(\beta)$-module $M$, the {\em $q$-character} of $M$ is defined by
$$ \ch_q(M) \seteq \sum_{\nu \in I^\beta} \dim_q (e(\nu)M) \nu ,$$
Here,  $ \dim_q V\seteq\sum_{k\in \Z} \dim(V_k)q^{k} $ for a graded vector space $V = \bigoplus_{k\in \Z} V_k $. It is well defined
whenever $\dim V_k<\infty$ for all $k\in\Z$.

For $i\in I$, let $L(\alpha_i)$ be the simple graded $R(\alpha_i)$-module
such that $\ch_qL(\al_i)=(i)$. For simplicity, we write $L(i)$ for $L(\alpha_i)$ if there is no confusion.

For graded $R(\beta)$-modules $M$ and $N $, let $\Hom_{R(\beta)}(M,N)$ be the
space of morphisms in $\Modg(R(\beta))$, i.e.\ the $\bR$-vector space of
homogeneous homomorphisms of degree $0$, and set
\addtolength{\my}{-10ex}
\eq&& \label{eq:tHom}
\hs{1ex}\tHom_{R(\beta)}(M,N)=\smash{\mathop\bigoplus\limits_{k \in \Z}} \tHom_{R(\beta)}(M,N)_k\\
&&\hs{35ex}\hfill\text{with $\tHom_{R(\beta)}(M,N)_k\seteq\Hom_{R(\beta)}(q^kM,N)$.}\nn
\eneq
We write $\tEnd_{R(\beta)}(M)$ for $\tHom_{R(\beta)}(M,M)$.
When $f \in \Hom_{R(\beta)}(q^{k}M, N)$, we denote
$$ \deg(f) := k. $$
For simplicity, we write $\tHom_R( M, N )$ for $\tHom_{R(\beta)}( M, N )$ if there is no confusion.

We write $[R\text{-proj}]$ and $[R\gmod]$ for
the (split) Grothendieck group of $R\proj$
and the Grothendieck group of $R\gmod$.
Then, the $\Z$-grading gives a $\Z[q,q^{-1}]$-module structure on
$[R\text{-proj}]$ and $[R\gmod]$,
and the convolution gives an algebra structure.

\Th [{\cite{KL09, KL11, R08}}]
There exist algebra isomorphisms
$$ [R\text{-{\rm proj}}] \simeq  U_\A^-(\g), \quad [R\gmod] \simeq  \Aq.$$
\enth
Here, $\Aq\seteq\set{a\in\Um}{(a,\Uam)\subset\A}$, where $(\cdot,\cdot)$ is the non-degenerate symmetric bilinear form
on $\Um$ defined in \cite{K91}. Note that $\Aq$ is an $\A$-subalgebra of
$\Uam$ (cf.\ \cite{KKKO18} where $\Aq$ is denoted by $\Aq[\mathfrak{n}]_{\Z[q^{\pm1}]}$).

\vskip 0.3em

\Def \label{dfn: grading}
Let $\fc$ be a $\Z$-valued skew-symmetric bilinear form on
$\rtl$.
If we redefine $\deg(\tau_l e(\nu)) $ to be
$-(\al_{\nu_l},\al_{\nu_{l+1}})-\fc(\al_{\nu_l},\al_{\nu_{l+1}})$,
then it gives a well-defined $\Z$-graded algebra structure on $R(\beta)$.
We denote by $R_\fc(\beta)$ the $\Z$-graded algebra thus defined.
\edf
The usual grading $\eqref{Eq: grading}$
is a special case of such  a $\Z$-grading.

We define $R_\fc(\beta)\gmod$, $R_\fc\gmod$, etc., similarly.

Let us denote by $\Modg(R_\fc(\beta))[q^{1/2}]$ the category of
 $(\frac{1}{2}\Z)$-graded modules over $R_\fc(\beta)$.
For $\nu\in I^\beta$ we set
$$H(\nu)=\frac{1}{2}\sum_{1\le a<b\le \height{\beta}}
\fc(\al_{\nu_a},\al_{\nu_b}).$$

\Lemma
For $\beta\in\rtlp$ and $M\in\Modg(R(\beta))[q^{1/2}]$,
set
$$\bl K_\fc(M)\br_n=\soplus_{\nu\in I^\beta}e(\nu)M_{n-H(\nu)}.$$
Then we have
\bnum
\item $K_\fc$ is an equivalence of categories
from $\Modg(R(\beta))[q^{1/2}]$ to $\Modg(R_\fc(\beta))[q^{1/2}]$.
\item For $M\in \Modg(R(\beta))[q^{1/2}]$ and $N\in\Modg(R(\gamma))[q^{1/2}]$,
we have
$$K_\fc(M\conv N)
\simeq q^{\frac{1}{2}\fc(\beta,\gamma)} K_\fc(M)\conv K_\fc(N).$$
\ee
\enlemma
\Proof
(i)\ We have
\begin{align*}
\tau_ke(\nu)\bl K_\fc(M)\br_n
&=\tau_ke(\nu)M_{n-H(\nu)} \\
& \subset e(s_k\nu)M_{n-H(\nu)-(\al_{\nu_k},\al_{\nu_{k+1}})}\\
&= e(s_k\nu)\bl K_\fc(M)\br_{n-H(\nu)-(\al_{\nu_k},\al_{\nu_{k+1}})+H(s_k\nu)}.
\end{align*}
Then (i) follows from
\begin{align*}
H(s_k\nu)-H(\nu)=\frac{1}{2}
\bl\fc(\al_{\nu_{k+1}},\al_{\nu_{k}})-\fc(\al_{\nu_{k}},\al_{\nu_{k+1}})\br
=-\fc(\al_{\nu_{k}},\al_{\nu_{k+1}}).
\end{align*}

\noi
(ii)\ For $\nu\in I^\beta$ and $\mu\in I^\gamma$, we have
\begin{align*}
e(\nu)K_\fc(M)_a\tens e(\mu)K_\fc(N)_b
&=e(\nu)M_{a-H(\nu)}\tens e(\mu)N_{b-H(\mu)}\\
&\subset e(\nu\ast\mu)(M\conv N)_{a+b-H(\nu)-H(\mu)}\\
&=e(\nu\ast\mu)K_\fc(M\conv N)_{a+b-H(\nu)-H(\mu)+H(\nu\ast\mu)}.
\end{align*}
Since
$$H(\nu\ast\mu)-H(\nu)-H(\mu)
=\frac{1}{2}\fc(\beta,\gamma),$$
we have
\eqn
K_\fc(M)_a\tens K_\fc(N)_b
\subset K_\fc(M\conv N)_{a+b+\frac{1}{2}\fc(\beta,\gamma)}.
\eneqn
Therefore we have
$$K_\fc(M)_a\tens K_\fc(N)_b\To
\bl q^{-\frac{1}{2}\fc(\beta,\gamma)}K_\fc(M\conv N)\br_{a+b},$$
which induces an isomorphism
$$K_\fc(M)\conv K_\fc(N)\isoto q^{-\frac{1}{2}\fc(\beta,\gamma)}K_\fc(M\conv N).$$
\QED

We define the algebra
$\Uam_\fc$ as $\Z[q^{\pm1/2}]\otimes_{\Z[q^{\pm1}]}\Uam$ with the new multiplication $\circ_\fc$
given by
$$a\circ_\fc b=q^{-\frac{1}{2}\fc(\al,\beta)}ab\quad
\text{for $a\in \Z[q^{\pm1/2}]\otimes_{\Z[q^{\pm1}]}\Uam_\al$ and
$b\in \Z[q^{\pm1/2}]\otimes_{\Z[q^{\pm1}]}\Uam_\beta$.}$$
We define $\Aq_\fc$ similarly.

\Cor\label{cor:twist}
There is a $\Z[q^{\pm1/2}]$-algebra isomorphism
$$\xi_\fc\col\Aq_\fc\isoto
\Z[q^{\pm1/2}]\tens_{\Z[q^{\pm1}]}[R_\fc\gmod].
$$
\encor

\subsection{Remark on parity}
Under hypothesis \eqref{item:inner}\,\eqref{hyp:even} in Definition \ref{def:cartan},
the category $\Modg(R(\beta))$ is divided
into two parts according to the parity of degrees for any $\beta\in\rtlp$.

\Lemma \label{lem:parity}
Let $\beta\in\rtlp$. Then there exists a map
$$S\cl I^\beta\To \Z/2\Z$$
such that
$$S(s_k\nu)=S(\nu)+(\al_{\nu_k},\al_{\nu_{k+1}})$$
for any
$\nu\in I^\beta$ and any integer $k$
with $1\le k<\Ht{\beta}$.
\enlemma
\Proof Set $n=\Ht{\beta}$.
Let us choose a total order $\prec$ on $I$.
Then we set
$$S(\nu)\seteq\sum_{1\le a<b\le n, \;\nu_{a}\prec\nu_b}(\al_{\nu_a},\al_{\nu_b}).$$
Then we have
\begin{align*}
S(s_k\nu)&=S(\nu)+\bl\delta(\nu_{k+1}\prec\nu_k)-
\delta(\nu_k\prec\nu_{k+1})\br(\al_{\nu_k},\al_{\nu_{k+1}})\\
&\equiv
S(\nu)+\bl1-\delta(\nu_{k}=\nu_{k+1})\br(\al_{\nu_k},\al_{\nu_{k+1}})
\equiv S(\nu)+(\al_{\nu_k},\al_{\nu_{k+1}})  \mod 2.
\end{align*}

Here, for a statement $P$, we set $\delta(P)$ to be $1$ if $P$ is true and $0$ if $P$ is false.
\QED

\Prop\label{prop:paritydec} Let $\beta\in\rtlp$ and $S\cl I^\beta\To \Z/2\Z$ be as in {\rm Lemma \ref{lem:parity}}.
Let $\Modg(R(\beta))^S$ be the full subcategory of
$\Modg(R(\beta))$ consisting of graded $R(\beta)$-modules $M$ such that
$e(\nu)M_k=0$ for any $\nu\in I^\beta$ and $k\equiv S(\nu)+1\;\mathrm{mod}\,2$.
Then we have
$$\Modg(R(\beta))\simeq \Modg(R(\beta))^S\soplus q\Modg(R(\beta))^S.$$
\enprop
\Proof
For any graded $R(\beta)$-module $M$ and $\eps=0,1$ set
$$M^\eps\seteq\soplus_{\substack{
\nu\in I^\beta,\;k\in\Z, \\ k\equiv S(\nu)+\eps\;\mathrm{mod}\,2}}e(\nu)M_k.
$$
Then, we can see easily that  they are $R(\beta)$-submodules of $M$
and $M=M^{0}\soplus M^{1}$.
 Moreover, we have $M^\eps\in q^{\eps}\Modg(R(\beta))^S$.
\QED

Note that $q^2\Modg(R(\beta))^S=\Modg(R(\beta))^S$ and
\eq
&&\text{$\tHom_{R(\beta)}(M,N)_k=0$ if $k$ is odd and  $M$, $N\in\Modg(R(\beta))^S$.}\label{eq:parity}
\eneq

\subsection{R-matrices}\label{Sec: R matrix}
Let $\beta\in\rtlp$ and $m=\Ht{\beta}$.
For $k = 1, \ldots, m-1$ and $\nu\in I^\beta$, the {\em intertwiner} $\varphi_k \in R(\beta)$ is defined by
$$ \varphi_k e(\nu) := \left\{
              \begin{array}{ll}
                (\tau_kx_k - x_k\tau_k) e(\nu) & \hbox{ if } \nu_k = \nu_{k+1}, \\
                \tau_k e(\nu) & \hbox{ otherwise.}
              \end{array}
            \right.
$$
\begin{lemma}[{\cite[Lem.\ 1.5]{K^3}}] \label{Lem: intertwiners} \
\bnum
\item $\varphi_k^2 e(\nu) = ( \qQ_{\nu_k, \nu_{k+1}} (x_k, x_{k+1}) + \delta_{\nu_k, \nu_{k+1}} )e(\nu). $
\item $\{ \varphi_k \}_{1 \le k \le m-1}$ satisfies the braid relation.
\item For a reduced expression $w = s_{i_1} \cdots s_{i_t}  \in \sg_{m}$, let $\varphi_w = \varphi_{i_1} \cdots \varphi_{i_t}$. Then $\varphi_w$ does not depend on the choice of reduced expressions of $w$.
\item For $w\in \sg_m$ and $1 \le k \le m$, we have $\varphi_w x_k = x_{w(k)} \varphi_w$.
\item For $w \in \sg_m$ and $1 \le k < m$, if $w(k+1) = w(k) + 1$, then $\varphi_w \tau_k = \tau_{w(k)} \varphi_w$.
\item $\varphi_{w^{-1}} \varphi_{w} e(\nu) = \prod_{a < b,\ w(a)>w(b)} ( \qQ_{\nu_a, \nu_{b}} (x_a, x_{b}) + \delta_{\nu_a, \nu_b} )e(\nu)$.
\ee
\end{lemma}

For $m,n \in \Z_{\ge0}$, let $w[m,n]$ be the element of $\sg_{m+n}$ defined by
$$
w[m,n](k) = \left\{
              \begin{array}{ll}
                k+n & \hbox{ if } 1 \le k \le m, \\
                k-m & \hbox{ if } m < k \le m+n.
              \end{array}
            \right.
$$
Let $M$ be an $R(\beta)$-module with $\Ht{\beta}=m$ and $N$ an $R(\beta')$-module with $ \Ht{\beta'}=n$.
The $R(\beta)\otimes R(\beta')$-linear map $M\otimes N \longrightarrow N \conv M$ given by $u\otimes v \mapsto \varphi_{w[n,m]}(v \otimes u)$ can be extended to the $R(\beta+\beta')$-module
homomorphism
\begin{align} \label{Eq: RMN}
R_{M,N}: M\conv N \longrightarrow N \conv M.
\end{align}

For $\beta = \sum_{k=1}^m \alpha_{i_k}$, we set $\supp(\beta)\seteq\{ i_k \mid 1 \le k \le m \}$.
\Def \label{def:symmetric}
The quiver Hecke algebra $R(\beta)$ is said to be {\em symmetric} if $\qQ_{i,j}(u,v)$ is a polynomial in $u-v$ for all $i,j \in \supp(\beta)$.
\edf

Suppose that $R(\beta)$ is symmetric. Let $z$ be an indeterminate. For an $R(\beta)$-module $M$, we define an $R(\beta)$-module structure on $M_z := \bR[z] \otimes_{\bR} M$ by
\begin{equation}\label{Eq: affinization1}
\begin{aligned}
&e(\nu) (a \otimes u) = a \otimes e(\nu)u, \quad
x_j(a \otimes u) = (za)\otimes u + a \otimes x_j u, \\
&\tau_k (a \otimes u) = a \otimes (\tau_k u),
\end{aligned}
\end{equation}
for $\nu \in I^n$, $a\in \bR[z]$ and $u\in M$.
We call $M_z$ the {\em affinization} of $M$.
For a non-zero $R(\beta)$-module $M$ and a non-zero $R(\beta')$-module $N$, let $s$ be the order of the zero of $R_{M_z, N_{z'}}: M_z \conv N_{z'} \rightarrow N_{z'} \conv M_z$, and
$$\Rm_{M_z, N_{z'}} := (z'-z)^{-s} R_{M_z, N_{z'}} .$$
We define $  \mathbf{r}_{M,N} : M \conv N \longrightarrow N\conv M  $ by
$$\mathbf{r}_{M,N} := \Rm_{M_z, N_{z'}}|_{z=z'=0}. $$

We put $R(\beta)[z_1, \ldots, z_k] := \bR[z_1, \ldots, z_k]\otimes_{\bR} R(\beta)$. For simplicity,
we write $R[z_1, \ldots, z_k]$ for $R(\beta)[z_1, \ldots, z_k]$ if there is no afraid of confusion.

\Th [{\cite[Section\,1]{K^3}}]  \label{Thm: R matrices}
Suppose that $R(\beta)$ and $R(\beta')$ are symmetric.
Let $M$ be a non-zero $R(\beta)$-module and $N$ a non-zero $R(\beta')$-module.
\bnum
\item $\Rm_{M_z, N_{z'}}$ and $\mathbf{r}_{M,N}$ are non-zero.
\item $\Rm_{M_z, N_{z'}}$ and $\mathbf{r}_{M,N}$ satisfy the braid relations.
\item Set
$$A = \sum_{\mu \in I^\beta, \nu \in I^{\beta'}}  \left(  \prod_{1 \le a \le m, 1 \le b \le n, \mu_a \ne \nu_b }
\qQ_{\mu_a, \nu_b}(x_a \boxtimes e(\beta'),e(\beta)\boxtimes x_b) \right)  e(\mu) \boxtimes e(\nu) .$$
Then, $A$ is in the center of $R(\beta)\otimes R(\beta')$, and
$$R_{N_{z'}, M_z} R_{M_z, N_{z'}} (u \otimes v)= A (u \otimes v) $$
for $u \in M_z$ and $v\in N_{z'}$.
\item If $M$ and $N$ are simple modules, then
\begin{align*}
\tEnd_{R(\beta + \beta')[z,z']} ( M_{z} \conv N_{z'} ) &\simeq \bR[z,z'],\\
\tHom_{R(\beta+\beta')[z,z']}( M_z \conv N_{z'}, N_{z'} \conv M_z ) &\simeq \bR[z,z'] \Rm_{M_z, N_{z'}}.
\end{align*}
\end{enumerate}
\enth

\section{Affinization}
\subsection{Definition of affinization}
\Def
For any $i\in I$ and $\beta\in \rtl_+$ with $\Ht{\beta}=m$, we set
$$\ap[i,\beta]=\sum_{\nu\in I^\beta}\Bigl(\prod_{a\in [1,m],\,\nu_a=i}x_a\Bigr)e(\nu),$$
where $[1,m] = \{1,2,\ldots,m \}$.
\edf
Note that $\ap[i,\beta]$ belongs to the center of $R(\beta)$.
If there is no afraid of confusion, we simply write
$\ap$ for  $\ap[i,\beta]$.

\Def\label{def:aff} Let $\beta\in\rtl_+$ and $\bM$  a simple $R(\beta)$-module.
An {\em affinization} $\Ma\seteq(\Ma,\z)$ of $\bM$ is an
$R(\beta)$-module $\Ma$ with an injective homogeneous endomorphism
$\z$ of\/ $\Ma$ of degree $\dM\in\Z_{>0}$
and an isomorphism $\Ma/\z\Ma\isoto \bM$
satisfying the following conditions\/{\rm:}

\bna
\item $\Ma$ is a finitely generated free module
over the polynomial ring $\cor[\z]$,\label{hyp:1}
\setcounter{myc}{\value{enumi}}
\item
$\ap\Ma\not=0$ for any $i\in I$.\label{hyp:nd}
\setcounter{mycounter}{\value{enumi}}
\ee
If it satisfies moreover the following additional condition,
we say that $\Ma$ is a {\em strong affinization}{\rm:}
\bna \setcounter{enumi}{\value{mycounter}}
\item the exact sequence
$0\to \z\Ma/\z^2\Ma\to\Ma/\z^2\Ma\to \Ma/\z\Ma\to0$ of $R(\beta)$-modules
does not split.\label{hyp:2}
\ee
We say that an affinization is {\em even} if $\dM$ is even.
\edf
Let us denote by $\pM\cl \Ma\epito\bM$ the composition
$\Ma\epito\Ma/\z\Ma\isoto\bM$.

\Rem\label{rem:aff}\hfill
\bnum
\item\label{it: affeq}
Condition \eqref{hyp:1} is equivalent to the condition

\setcounter{enumi}{\value{myc}}
\be[{(\rm a')}]
%(\ref{hyp:1}')
\item \label{hyp:1'}The degree of $\Ma$ is bounded from below, that is,
$\Ma_n=0$ for $n\ll0$.
\ee
Moreover, under these equivalent conditions, we have
$$\ch_q(\Ma)=(1-q^{\dM})^{-1}\ch_q(\bM).$$
Note that a finitely generated $R$-module $M$ satisfies the condition
(\ref{hyp:1}').
\item
The non-splitness condition \eqref{hyp:2} is equivalent to saying that
$\z\Ma/\z^2\Ma$ is a unique proper $R(\beta)$-submodule of $\Ma/\z^2\Ma$.
\item If $R(\beta)$ is a symmetric quiver Hecke algebra, then $\bM_z$ is a strong affinization of any simple $R(\beta)$-module $\bM$ for $\beta\not=0$.
\ee
\enrem

\begin{example}
\bnum
\item For $i\in I$,
$\Ma\seteq L(i)_z\conv L(i)$ is not an affinization of $\bM\seteq L(i)\conv L(i)$
In fact, conditions \eqref{hyp:1} and \eqref{hyp:2}
in Definition~\ref{def:aff} hold
but condition \eqref{hyp:nd} does not.

\item
Let $(\Ma,\z)$ be an affinization of $\bM$.
Assume that $\dM=ab$ for $a,b \in \Z_{>0}$
and let $z$ be an indeterminate of homogeneous degree $b$.
Let $\cor[\z]\to\cor[z]$ be the algebra homomorphism given by $\z\mapsto z^a$.
Then $\bl\cor[z]\tens_{\cor[\z]}\Ma,\,z\br$
is an affinization of $\bM$.
If $a>1$ then it is never a strong affinization,
because $$\bl\cor[z]\tens_{\cor[\z]}\Ma\br/\bl z^a\cor[z]\tens_{\cor[\z]}\Ma\br
\simeq (\cor[z]/\cor[z]z^a)\tens_\cor\bM$$
is a semisimple $R(\beta)$-module.
\ee
\end{example}

As seen in the proposition below, every affinization is essentially even.

\Prop
Let $(\Ma,\z)$ be an affinization of a simple module $\bM$.
Assume that the homogeneous degree $\dM$ of $\z$ is odd.
Then there exists an $R(\beta)$-submodule $\Ma'$ of $\Ma$ such that
\bnum
\item $\z^2\Ma'\subset \Ma'$, and
$(\Ma',\,\z^2)$ is an affinization of $\bM$,
\item $\Ma\simeq \cor[\z]\tens_{\cor[\z^2]}\Ma'$ as an $\R(\beta)[\z]$-module.
\ee
\enprop
\Proof
Let $\Modg(R(\beta))\simeq \Modg(R(\beta))^S\soplus q\Modg(R(\beta))^S$
be the decomposition in Proposition~\ref{prop:paritydec}.
We may assume that $\bM$ belongs to $\Modg(R(\beta))^S$.
Let $\Ma=\Ma'\soplus \Ma''$ be the decomposition with
$\Ma'\in \Modg(R(\beta))^S$ and $\Ma''\in q\Modg(R(\beta))^S$.
Then $\z \Ma'\subset \Ma''$ and $\z \Ma''\subset\Ma'$ by \eqref{eq:parity}.
Hence, we have the decomposition
$$  \Ma /\z \Ma=(\Ma'/\z\Ma'')\soplus (\Ma''/\z\Ma'),$$
which implies that
$\Ma'/\z\Ma''\simeq\bM$ and $\Ma''=\z\Ma'$.
Thus we obtain the desired result.
\QED

\subsection{Strong affinization}
Note that Lemmas~\ref{le:affstrong} and \ref{lem:affweek} below
hold without assumption \eqref{hyp:nd} in Definition~\ref{def:aff}.

\Lemma\label{le:affstrong}
Assume that
\begin{equation*}
\addtocounter{equation}{1}
\setcounter{mycounter}{\value{section}}
\setcounter{mycounters}{\value{equation}}
\hs{-3ex}(\themycounter.\themycounters)_{\rm strong}\hs{1ex}
\left\{\parbox{70ex}{$\beta\in\rtl_+$ and $(\Ma,\z)$ is a {\em strong} affinization of a simple
$R(\beta)$-module $\bM$, $\z$ has homogeneous degree $\dM\in\Z_{>0}$,
and $\pM\cl \Ma\to \bM$ is a canonical projection.}\right.
\end{equation*}
\bnum
\item The head of the $R$-module $\Ma$ is isomorphic to $\bM$,
or equivalently $\z\Ma$ is a unique maximal $R(\beta)$-module. \label{le:afhead}
\item Let $s\seteq\min\set{m\in\Z}{\Ma_m\not=0}$, and
$u$ a non-zero element of $\Ma_s$. Then, $\Ma=R(\beta)u$. \label{lem:cyc}
\item  $\tEnd_{R(\beta)}(\Ma)\simeq\cor[\z]\id_{\Ma}$. \label{lem:staffend}
\ee
\enlemma
\Proof
(i) Let $S$ be a simple module and $\vphi\cl \Ma \to S$ be an epimorphism.
By consideration of the homogeneous degree, we may assume that
$\vphi(\z^k\Ma)=0$ for $k\gg0$.
Let us take $k\ge0$ such that $\vphi(\z^k\Ma)=S$ and $\vphi(\z^{k+1}\Ma)=0$.
Since $\z^k\Ma/\z^{k+1}\Ma \simeq \Ma/\z\Ma$ is simple,
$\vphi$ induces an isomorphism $\z^k\Ma/\z^{k+1}\Ma\isoto S$.
It is enough to show that $k=0$.
If $k>0$ then we have a commutative diagram
\eq
&&\xymatrix
{0\ar[r]& \z^{k}\Ma/\z^{k+1}\Ma\ar[r]\ar[dr]^-\sim& \z^{k-1}\Ma/\z^{k+1}\Ma\ar[r]
\ar[d]^\vphi&\z^{k-1}\Ma/\z^{k}\Ma
\ar[r]& 0\\
&&S} .
\eneq
Hence the first row of the above diagram is a split exact sequence,
which contradicts Definition~\ref{def:aff}~\eqref{hyp:2}.

(ii) Since $u\not\in \z\Ma$, (i) implies that $\Ma=R(\beta)u$.

(iii) Let $f\in\tEnd_{R(\beta)}(\Ma)$ be a homogeneous endomorphism of degree $\ell$.
Assume that $f(\Ma)\subset \z^k\Ma$ for $k\in\Z_{\ge0}$.
We shall show $f\in \cor[\z]\id_{\Ma}$ by the descending induction on $k$.
If $\dM k>\ell$, then $f$ should be $0$ since $f(u)\not\in \z^k\Ma$
if $f(u)\not=0$. Here $u$ is as in (ii).
Suppose that $\dM k \le \ell$.
As $\bM$ is the head of $\Ma$,
the composition $\Ma\To[\z^{-k}f]\Ma\To[\pM] \bM$
decomposes as $\Ma\To[\pM]\bM\to\bM$. Hence the composition
must be equal to $c\pM$ for some $c\in \cor$,
which yields
$$ ( \z^{-k}f - c \id_{\Ma} ) (\Ma) \subset \z \Ma. $$
Therefore, we have $(f-c\z^k)(\Ma)\subset \z^{k+1}\Ma$,
and the induction proceeds.
\QED

\subsection{Normalized R-matrices}

\Lemma\label{lem:affweek}
Assume that
\begin{equation*}
\addtocounter{equation}{1}
\setcounter{mycountert}{\value{section}}
\setcounter{mycounterf}{\value{equation}}
\hs{-3.5ex}
(\themycountert.\themycounterf)_{\rm weak}\hs{1ex}\left\{
\parbox{70ex}{$\beta\in\rtl^+$ and $(\Ma,\z)$ is an affinization of a simple
$R(\beta)$-module $\bM$, $\z$ has homogeneous degree $\dM\in\Z_{>0}$,
and $\pM\cl \Ma\to \bM$ is a canonical projection.}\right.
\end{equation*}
\bnum
\item \label{lem:affend} $\tEnd_{R(\beta)[\z]}(\Ma)\simeq\cor[\z]\id_{\Ma}$.
\item \label{cor:MNNM}  For any $i\in I$, there exist $c_i\in\cor^\times$ and $d_i\in\Z_{\ge0}$ such that
$\ap\vert_\Ma=c_i\z^{d_i}$.
\ee
\enlemma
\Proof
(i) The proof is similar to that of Lemma \ref{le:affstrong} $\eqref{lem:staffend}$.
Let $f\in\tEnd_{R(\beta)[\z]}(\Ma)$ be a homogeneous endomorphism of degree $\ell$.
Suppose that $f(\Ma)\subset \z^k\Ma$ for $k\in\Z_{\ge0}$. We shall show $f\in \cor[\z]\id_{\Ma}$ by the descending induction on $k$.

We have $f=0$ if $\dM k>\ell$ by the degree consideration.
If $\dM k \le\ell$, then
the endomorphism $\z^{-k}f$ induces an endomorphism of $\bM$.
Hence it must be equal to $c\id_{\bM}$ for some $c\in \cor$.
Then $(f-c\z^k)(\Ma)\subset \z^{k+1}\Ma$,
and the induction proceeds.

(ii) The assertion follows from $\eqref{lem:affend}$ immediately.
\QED

\Lemma\label{lem:fg}
Let $\beta$, $\Ma$ and $\bM$ be as in
$(\themycountert.\themycounterf)_{\rm weak}$.
Assume further that $\beta\not=0$.
Then $\Ma$ is a finitely generated $\R(\beta)$-module.
\enlemma
\Proof
Since $\beta\not=0$, there exists $i\in I$ such that
$\ap[i,\beta]$ has a positive degree.
Then there exists $m>0$ such that $\z^m\in\cor\, (\ap[i,\beta]\vert_\Ma)\subset
\tEnd_R(\Ma)$.
Since $\Ma$ is finitely generated over $\cor[\z^m]$,
we obtain the desired result.
\QED

\Lemma\label{lem:Riso}
Let $\beta$, $\Ma$ and $\bM$ be as in
$(\themycountert.\themycounterf)_{\rm weak}$.
Let $\gamma\in\rtl_+$ and $N\in R(\gamma)\gmod$.
\bnum
\item
The homomorphisms
$$\text{$R_{\Ma[\z^{-1}],N}\cl\Ma [\z^{-1}] \conv N\to N\conv \Ma[\z^{-1}]$ and
$R_{N,\Ma[\z^{-1}]}\cl N\conv \Ma[\z^{-1}]\to\Ma [\z^{-1}]\conv N$}$$
are isomorphisms.
Here, $\Ma[\z^{-1}] = \bR[\z, \z^{-1}]\otimes_{\bR[\z]}\Ma$.
\item If $N$ is a simple module, there exist $c\in\cor^\times$ and $d\in\Z_{\ge0}$ such that
$R_{N,\Ma}\conv R_{\Ma,N}=c(\z^d\conv N)$ and
$R_{\Ma,N}\conv R_{N,\Ma}=c(N\conv \z^d)$.  \label{lem:comR}
\ee
\enlemma
\Proof
(i) is an immediate consequence of (ii).
Let us show (ii).
Set $m=\Ht{\beta}$ and $n=\Ht{\gamma}$.
Then, $(R_{N,\Ma}\conv R_{\Ma,N})\vert_{\Ma\tens N}$ is given by
$$\sum_{\nu\in I^{\beta+\gamma}}\Bigl(\prod_{1\le a\le m<b\le m+n, \ \nu_a \ne \nu_b } \qQ_{\nu_a,\nu_b}(x_a,x_b)\Bigr)e(\nu).$$
Since any element in the center of $R(\gamma)$ with positive degree acts by zero on $N$,
it is equal to
$$\sum_{\nu\in I^{\beta+\gamma}}\Bigl(\prod_{1\le a\le m<b\le m+n , \ \nu_a \ne \nu_b }\qQ_{\nu_a,\nu_b}(x_a,0)\Bigr)e(\nu).$$
Hence it is equal to a product of $\ap[i,\beta]\vert_\Ma$'s
up to a constant multiple.
Hence Lemma~\ref{lem:affweek} \eqref{cor:MNNM} implies the desired result.
\QED

Let $\Ma$ and $\bM$ be as in $(\themycountert.\themycounterf)_{\rm weak}$,
and let $N\in R\gmod$ be a non-zero
module.
Let $s$ be the largest integer such that
$R_{\Ma,N}(\Ma\conv N)\subset \z^sN\conv \Ma$.
Then we set
$$
\Rm_{\Ma,N} =\z^{-s}R_{\Ma,N}\cl \Ma\conv N\to N\conv \Ma.
$$
We denote by
$$
\rmat{\bM,N}\cl \bM\conv N\to N\conv \bM
$$
the homomorphism induced by $\Rm_{\Ma,N}$.
By the definition, $\rmat{\bM,N}$ never vanishes.
We set $\Rm_{\Ma,N}=0$, and $\rmat{\bM,N}=0$ when $N=0$.

We define $\Rm_{N,\Ma}$ and $\rmat{N,\bM}$ similarly.

\medskip
The arguments in \cite{KKKO18,KKKO14} still work well under these assumptions,
and we obtain similar results.
We list some of them without repeating the proof.
A simple module $S$ is called {\em real} if $S\conv S$ is simple.

\Prop[\protect{\cite[Th.\ 3.2, Prop.\ 3.8]{KKKO14}, \cite[Prop.\ 3.2.9, Th.\ 4.1.1]{KKKO18}}] \label{Prop: simple-unique}
Assume that
\eq
&&\left\{\parbox{75ex}{
\bna
\item
$M$ and $N$  are simple $R$-modules,
\item one of them is real simple and also admits an affinization.
\ee}\right.\label{hyp:affMN}
\eneq
\bnum
\item $M\conv N$ has a simple head and a simple socle.
Moreover, $\Im(\rmat{M,N})$ is equal to the head of $M\conv N$ and the socle of $N\conv M$. \label{prop:simple}
\item We have
$$\tHom_R(M\conv N,M\conv N)=\cor\id_{M\circ N}$$
and
$$\tHom_R(M\conv N,N\conv M)=\cor\,\rmat{M,N}.$$  \label{prop:unique}
\item $M\hconv N$ appears only once in a
Jordan-H\"older series of $M\conv N$ in $R\smod$.
\ee
\enprop

\Prop \label{prop:uniqueR}
Let $\Ma$ and $\bM$ be as in $(\themycountert.\themycounterf)_{\rm weak}$,
and let $N$ be a simple $R$-module.
Assume that $\bM$ is real. Then we have
\bnum
\item
\eq
&&\tHom_{R[\z]}(\Ma\conv N,\Ma\conv N)=\cor[\z]\id_{\Ma\circ N}
\label{eq:uniend1}
\eneq
and
\eq&&\tHom_{R[\z]}(N\conv\Ma, N\conv \Ma)=\cor[\z]\id_{N\circ\Ma}.
\label{eq:uniend2}
\eneq
\item \hskip 0.05em  $\tHom_{R[\z]}(\Ma\conv N, N\conv\Ma)$ and $\tHom_{R[\z]}(N\conv\Ma,  \Ma\conv N)$ are free $\cor[\z]$-modules of rank one.
\ee
\enprop
\Proof
(i) Let us first show \eqref{eq:uniend1}.
The idea for the proof is similar to that of Lemma \ref{le:affstrong} $\eqref{lem:staffend}$.

Let $f\in\tHom_{R[\z]}(\Ma\conv N,\Ma\conv N)$ of homogeneous degree $\ell$.
We know that $f(\Ma\conv N)\subset \z^s\Ma\conv N$ for some $s\in \Z_{\ge0}$.
We shall show $f\in \cor[\z]\id_{\Ma\circ N}$ by the descending induction on $s$.
If $s\gg0$, then $f$ should be zero by the degree consideration.
Now, we consider $\z^{-s}f$.
As $\z^{-s}f$ induces an endomorphism of $\bM\conv N$,
by Proposition \ref{Prop: simple-unique} $\eqref{prop:unique}$, it is equal to $c\id_{\bM\circ N}$ for $c \in \bR$.
Hence we have
$$
\bl f-c z^s_{\Ma}\br(\Ma\conv N)\subset \z^{s+1}\Ma\conv N.
$$
Thus, the induction hypothesis implies that
$f-c z^s_{\Ma}\in \cor[\z]\id_{\bM\circ N}$.
The proof of \eqref{eq:uniend2} is similar.

(ii) By Lemma~\ref{lem:Riso}, we have an $R[\z]$-linear monomorphism
$N\conv\Ma\monoto \Ma\conv N$.
Hence
we have
$$\tHom_{R[\z]}(\Ma\conv N,N\conv\Ma)\monoto
\tHom_{R[\z]}(\Ma\conv N,\Ma\conv N)\simeq \cor[\z].$$
As $\tHom_{R[\z]}(\Ma\conv N,N\conv\Ma)$ is non-zero,
$\tHom_{R[\z]}(\Ma\conv N,N\conv\Ma)$ is a free $\cor[\z]$-module of rank one.
\QED

\Prop \label{prop:MNNM}
We assume that
\eq
&&\left\{\parbox{70ex}
{\bna
\item $(\Ma,\z)$ and $(\Na,\zN)$ are affinizations
of simple modules $\bM$ and $\bN$, respectively,
\item one of $\bM$ and $\bN$ is real.
\ee}\right.\label{hyp:MN}
\eneq
Then we have
\bnum
\item  $\tHom_{R[\z,\zN]}(\Ma\conv \Na,\Ma\conv \Na) =\cor[\z,\zN]\id_{\Ma\circ \Na}$,
\item $\tHom_{R[\z,\zN]}(\Ma\conv \Na,\Na\conv \Ma)$ is a free $\cor[\z,\zN]$-module of rank one.
\ee
\enprop
\Proof
(i) Assume that $\bM$ is real simple. The other case can be proved similarly.

Let $f$ be a homogeneous element of
$\tHom_{R[\z,\zN]}(\Ma\conv \Na,\Ma\conv \Na)$ of degree $\ell$. Assuming that $\Im(f)\subset \zN^k(\Ma\conv \Na)$, we shall show
$f\in \cor[\z,\zN]\id_{\Ma\circ \Na}$ by the descending induction on $k$.
If $k\gg0$, then $f$ should be zero by the homogeneous-degree consideration.
We now consider $\zN^{-k}f$. The $R[\z,\zN]$-linear homomorphism
$\zN^{-k}f\cl \Ma\conv \Na\to\Ma\conv \Na$ induces an
$R[\z]$-linear homomorphism
$\Ma\conv \bN\to\Ma\conv\bN$.
By Proposition~\ref{prop:uniqueR},
it is equal to $\vphi(\z)\id_{\Ma\circ\bN}$ for some $\vphi(\z)\in\cor[\z]$.
Hence we have
$$\Im(f-\zN^k\vphi(\z)\id_{\Ma\circ\Na})\subset \zN^{k+1}\Ma\conv \Na,$$
which implies
$$f-\zN^k\vphi(\z)\id_{\Ma\circ \Na}\in\cor[\z,\zN]\id_{\Ma\circ \Na}$$
by the induction hypothesis.

\smallskip\noi 
(ii)
Since $R_{\Ma,\Na} R_{\Na,\Ma} |_{\zN = 0}$ is non-zero by Lemma \ref{lem:Riso} $\eqref{lem:comR}$, the assertion (i) tells that
$R_{\Ma,\Na} R_{\Na,\Ma} \in \cor[\z,\zN]\id_{\Na\circ \Ma}$ is non-zero.
The injectivity of $R_{\Ma,\Na} R_{\Na,\Ma}$ implies that $R_{\Na,\Ma}\cl\Na\conv \Ma\monoto\Ma\conv\Na$ is injective. Thus the composition with $R_{\Na,\Ma}$ induces
an injective homomorphism
$$\tHom_{R[\z,\zN]}(\Ma\conv \Na,\Na\conv \Ma)
\monoto \tHom_{R[\z,\zN]}(\Ma\conv \Na,\Ma\conv \Na)
\simeq\cor[\z,\zN].$$

We now consider the non-zero $\cor[\z,\zN]$-module $L\seteq\tHom_{R[\z,\zN]}(\Ma\conv \Na,\Na\conv \Ma)$.
Let $a,b\in \cor[\z,\zN]$ be non-zero relatively prime elements. If $f \in aL \cap bL$, then we have
$ f(\Ma \circ \Na) \subset a(\Na \circ \Ma) \cap b(\Na \circ \Ma)$.
Since $a$ and $b$ are relatively prime and $\Na \circ \Ma$ is a free
$\cor[\z,\zN]$-module,
$$ f(\Ma \circ \Na) \subset ab(\Na \circ \Ma),$$
which implies
that $(ab)^{-1} f : \Ma \circ \Na \rightarrow \Na \circ \Ma$ is well-defined, i.e., $f \in ab L$.
Therefore, we conclude that $L$ satisfies the condition:
\eqn&
\parbox{75ex}{if non-zero elements $a,b\in \cor[\z,\zN]$
are prime to each other,
then we have
$aL\cap bL=abL$,}
\eneqn
which yields that $L$ is a free $\cor[\z,\zN]$-module of rank one.
\QED

We define
$\Rm_{\Ma,\Na}$ as a generator of
the $\cor[\z,\zN]$-module $\tHom_{R[\z,\zN]}(\Ma\conv \Na,\Na\conv \Ma)$.
It is uniquely determined up to a constant multiple.
We call it a {\em normalized R-matrix}.

\Th \label{Thm: specializetion of Rm}
Assume \eqref{hyp:MN}.
Then,
$\Rm_{\Ma,\Na}\vert_{\z=\zN=0}\cl\bM\conv\bN\to\bN\conv\bM$
 does not vanish and is equal to $\rmat{\bM,\bN} $
up to a constant multiple.
\enth
\Proof
Since any simple $R$-module is absolutely simple,
we may assume that the base field $\cor$ is algebraically closed
without loss of generality.

By Proposition~\ref{Prop: simple-unique} $\eqref{prop:unique}$, we have
\eq
\tHom_{R}(\bM\conv \bN,\bN\conv \bM)=\cor\,\rmat{\bM,\bN}
\label{eq:rmat}
\eneq
for a non-zero $\rmat{\bM,\bN}\in \tHom_{R}(\bM\conv \bN,\bN\conv \bM)$.
Let $\ell$ be the homogeneous degree of $\rmat{\bM,\bN}$.

For $a\in\Z$, let $\cor[\z,\zN]_a$ be the homogeneous part
of $\cor[\z,\zN]$ of degree $a$ and set
$\cor[\z,\zN]_{\ge a}=\soplus_{k\ge a}\cor[\z,\zN]_k$.
Let $c\in\Z$ be the largest integer such that
$$\Rm_{\Ma,\Na}(\Ma\conv \Na)\subset \cor[\z,\zN]_{\ge c}(\Na\conv \Ma).$$
Then $\Rm_{\Ma,\Na}$ induces a non-zero map
$$\vphi\cl \bM\conv\bN\To \Bigl(\cor[\z,\zN]_{\ge c}(\Na\conv \Ma)\Bigr)/
\Bigl(\cor[\z,\zN]_{\ge c+1}(\Na\conv \Ma)\Bigr).$$
Since $\Na\conv \Ma$ is a free $\cor[\z,\zN]$-module, we have
$$\Bigl(\cor[\z,\zN]_{\ge c}(\Na\conv \Ma)\Bigr)/
\Bigl(\cor[\z,\zN]_{\ge c+1}(\Na\conv \Ma)\Bigr)
\simeq \cor[\z,\zN]_{c}\tens(\bN\conv \bM).$$
By \eqref{eq:rmat}, there exists a non-zero $f(\z,\zN)\in \cor[\z,\zN]_{c}$
such that
$$\text{$\vphi(u)=f(\z,\zN)\rmat{\bM,\bN}(u)$ for any $u\in\bM\conv\bN$.}$$
Hence, the homogeneous degree of $\Rm_{\Ma,\Na}$ is equal to
$c+\ell$, and
$$\Rm_{\Ma,\Na}(\Ma\conv \Na)\subset f(\z,\zN)(\Na\conv \Ma)+\cor[\z,\zN]_{\ge c+1}(\Na\conv \Ma).$$
Let us show that $f(\z,\zN)$ is a constant function (i.e., $c=0$).
Assuming that $c>0$, let us take
a prime divisor $a(\z,\zN)$ of $f(\z,\zN)$.
Let $(x,y)$ be an arbitrary point of $\cor^2$
such that $a(x,y)=0$.
Let $\dM$ and $\dM[\Na]$ be the homogeneous degree of $\z$ and $\zN$,
respectively,
and let $z$ be an indeterminate of homogeneous degree one.
Then $a(xz^{\dM},yz^{\dM[\Na]})=f(xz^{\dM},yz^{\dM[\Na]})=0$.
Let $\cor[\z,\zN]\to\cor[z]$ be the map obtained by the substitution
$\z=xz^{\dM}$ and $\zN=yz^{\dM[\Na]}$.
Set
\begin{align*}
K=\cor[z]\tens_{\cor[\z,\zN]}(\Ma\conv\Na),\quad
K'=\cor[z]\tens_{\cor[\z,\zN]}(\Na\conv\Ma), \quad R'=\cor[z]\tens_{\cor[\z,\zN]}\Rm_{\Ma,\Na}.
\end{align*}
Then we obtain the map
$$R'\cl K\to z^{c+1}K'.$$
Note that $K/zK\simeq\bM\conv \bN$ and $K'/zK'\simeq\bN\conv\bM$.
We shall show $R'(K)\subset z^{k}K'$
for any $k\ge c+1$ by induction on $k$.
Assume that $k\ge c+1$ and
$R'(K)\subset z^{k}K'$.
Then the morphism $ \bM \conv\bN \rightarrow \bN \conv \bM $ induced by $z^{-k}R'$ is equal to $b\,\rmat{\bM,\bN}$ for some $b \in \cor$.
If $b\not=0$, then the homogeneous degree of $R'$
is equal to $k+\ell>c+\ell$, which is a contradiction.
Thus $b=0$ and
$R'(K)\subset z^{k+1}K'$.
Hence the induction proceeds, and we conclude that
$$\Rm_{\Ma,\Na}\vert_{\z=xz^{\dM},\,\zN=yz^{\dM[\Na]}}=0$$
for any $(x,y) \in \bR^2$ such that $a(x,y)=0$,
which implies $\Rm_{\Ma,\Na}$ is divisible by $a(\z,\zN)$.
This is a contradiction.

Therefore $f$ is a constant function, and
$\Rm_{\Ma,\Na}$ induces $\rmat{\bM,\bN}$ (up to a constant multiple)
after the specialization $\z=\zN=0$.
\QED

\Cor\label{cor:Rnorm} Assume \eqref{hyp:MN}.
If $\bM$ is real, then $\Rm_{\Ma,\Na}\vert_{\zN=0}=\Rm_{\Ma,\bN}$
\ro up to a constant multiple\rf.
\encor

\Lemma\label{lem:monic}
Assume \eqref{hyp:MN}. Then there exists a homogeneous element $f(\z,\zN)$
such that
\bnum
\item
$\Rm_{\Na,\Ma}\circ\Rm_{\Ma,\Na}=f(\z,\zN)\id_{\Ma\circ\Na}$ and
$\Rm_{\Ma,\Na}\circ\Rm_{\Na,\Ma}=f(\z,\zN)\id_{\Na\circ\Ma}$,
\item $f(\z,0)$ and
$f(0,\zN)$ are non-zero.
\ee
\enlemma
\Proof
It follows from Proposition~\ref{prop:MNNM},
Corollary~\ref{cor:Rnorm} and Lemma~\ref{lem:Riso}.
\QED

\Lemma \label{Lem: affMN}
Let $(\Ma,\z)$ and $(\Na,\zN)$
be affinizations of simple modules $\bM$ and $\bN$, respectively.
Assume that either $\bM$ or $\bN$ is real, and
$\bM\conv\bN\simeq\bN\conv \bM$.
Let $d$ be a common divisor of the homogeneous degree $\dM$ of $\z$
and $\dN$ of $\zN$.
Let $z$ be an indeterminate of homogeneous degree $d$ and let
$\cor[\z,\zN]\to\cor[z]$ be the algebra homomorphism given by
$\z\mapsto z^{\dM/d}$ and $\zN\mapsto z^{\dN/d}$.
Then, $\Ma\underset{z}{\conv}\Na\seteq\cor[z]\tens_{\cor[\z,\zN]}(\Ma\conv \Na)$ 
is an affinization of $\bM\conv\bN$.
\enlemma
\Proof By the condition, $\bM\conv\bN$ is simple.
Condition \eqref{hyp:1} in Definition~\ref{def:aff}
is obvious. Condition~\eqref{hyp:nd} follows from
$\ap\vert_{\Ma\underset{z}{\circ}\Na}=(\ap\vert_\Ma)\underset{z}{\conv}(\ap\vert_\Na)$. 
\QED

\Prop\label{prop:MM}
Let $\Ma$ and $\bM$ be as in $(\themycountert.\themycounterf)_{\rm weak}$.
Assume that $\bM$ is real.
We normalize $\Rm_{\Ma,\Ma}$ so that
it induces $\id_{\bM\circ\bM}$ after the specialization
$\z \conv \Ma= \Ma\conv\z=0$.
Then, we have
\bnum
\item
$(\Rm_{\Ma,\Ma}-\id_{\Ma\circ\Ma})(\Ma\conv \Ma)\subset
\bl\z \conv \Ma - \Ma\conv\z\br(\Ma\conv \Ma).$
\item
$\Rm_{\Ma,\Ma}\circ\Rm_{\Ma,\Ma}=\id_{\Ma,\Ma}$.
\ee
\enprop
\Proof
(i) To avoid possible confusions, let $(\Na,\zN)$ be a copy of $(\Ma,\z)$
and we regard $\Rm_{\Ma,\Ma}$ as a homomorphism $\Ma\conv \Na \rightarrow \Na \conv \Ma$.
We denote by $\imath \cl\Ma\conv \Na\isoto\Na\conv \Ma$ the identity.
We regard $\Ma\conv \Na$ and $\Na\conv \Ma$ as $R[\z,\zN]$-modules.
Then, $\Rm_{\Ma,\Na}$ commutes with $\z$ and $\zN$, but $\imath$ does not.
Precisely, we have $\imath\circ \z=\zN\circ \imath$ and $\imath\circ \zN=\z\circ \imath$.

Let $z$ be another indeterminate with the same homogeneous degree $\dM$,
and let
$$\cor[\z,\zN]\to\cor[z]$$
be the algebra homomorphism
give by $\z\mapsto z$ and $\zN\mapsto z$.
Then, Lemma \ref{Lem: affMN} implies that
$K\seteq \cor[z]\tens_{\cor[\z,\zN]}(\Ma\conv \Na)$ is an affinization of
$\bM\conv\bM$.
The homomorphisms $\Rm_{\Ma,\Na}$ and $\imath$ induce $R[z]$-linear endomorphisms
$R'$ and $\imath'$ of
$K$. By Lemma~\ref{lem:affweek} $\eqref{lem:affend}$, $R'$ and $\imath'$ are powers of $z$
up to a constant multiple.
Since they are $\id_{\bM\circ\bM}$ after the specialization $z=0$, we conclude
that $R'=\imath'$, which completes the proof.

(ii) The assertion follows from Lemma~\ref{lem:monic} immediately.
\QED

\Ex\label{Ex:Affn}
Let $i\in I$.
Let $P(i^n)$ be a projective cover of the simple module $ L(i)^{\circ\, n}$.
Then $P(i^n)$ is an $R(n\al_i)$-module generated by an element $u$ 
of degree $0$ 
with the defining relation
$\tau_ku=0$ ($1\le k<n$).
Let $e_k(x_1,\ldots, x_n)$ be the elementary symmetric function of degree $k$.
The center of $R(n\al_i)$ is equal to
$\cor[e_k(x_1,\ldots, x_n)\mid k=1,\ldots,n]=\cor[x_1,\ldots,x_n]^{\sym_n}$.
Then we have
$$L(i)^{\circ\,n}\simeq P(i^n)/\Bigl(\sum_{k=1}^{n}R(n\al_i)e_k(x_1,\ldots, x_n)u\Bigr).$$
Set
$$\KM(i^n)\seteq P(i^n)/\Bigl(\sum_{k=1}^{n-1}R(n\al_i)e_k(x_1,\ldots, x_n)u\Bigr),$$
and define $\zj[{\KM(i^n)}]\in\tEnd_{R(n\al_i)}(\KM(i^n))$ by
$$\zj[{\KM(i^n)}]u=e_n(x_1,\ldots, x_n)u.$$
Then $(\KM(i^n), \zj[{\KM(i^n)}])$ is a strong affinization of $L(i)^{\circ\, n}$.
Note that $\ap\vert_{\KM(i^n)}=\zj[{\KM(i^n)}]$.
The homogeneous degree of $\zj[{\KM(i^n)}]$ is $n(\al_i,\al_i)$.
\enEx
\section{Root modules}\label{Sec: root modules}\
In this section, we shall review results of
McNamara (\cite{Mc12}) and Brundan-Kleshchev-McNamara (\cite{BKM14}).
Throughout this section, we assume that the Cartan matrix $\cmA$ is of finite type.
Fix a reduced expression $w_0 = r_{i_1}r_{i_2} \ldots r_{i_N}$ of the longest element $w_0 \in \weyl$.
This expression gives a convex total order $\prec$ on the set
$\pr$ of positive roots: $\al_{i_1}\prec r_{i_1}\al_{i_2}\prec\cdots\prec
r_{i_1}\cdots r_{i_{N-1}}\al_{i_N}$.
 To each positive root $\beta \in \pr$, McNamara defined
a simple $R(\beta)$-module $L(\beta)$,
which he called the {\em cuspidal module} (\cite{KR11,Mc12}).

\Lemma[{\cite[Lem.\ 3.4]{Mc12}}] For any $\beta\in \pr $,
$L(\beta)$ is a real simple module.
\enlemma
\begin{lemma} [{\cite[Lem.\ 3.2]{BKM14}}] For $n\ge 0$, there exist unique
\ro up to isomorphism\rf\  $R(\beta)$-modules $\Delta_n(\beta)$ with $\Delta_0(\beta) = 0$ such that
there are short exact sequences
\begin{align*}
&0 \To q_\beta^{2(n-1)} L(\beta) \To[\ i_n\ ] \Delta_n(\beta) \To[\ p_n\ ] \Delta_{n-1}(\beta) \To 0,\\
&0 \To q_\beta^{2} \Delta_{n-1}(\beta) \To[\ j_n\ ]\Delta_n(\beta) \To[\ q_n\ ]L(\beta)
\To 0 \hs{6ex}\text{for $n\ge1$,}
\end{align*}
where $q_\beta = q^{(\beta,\beta)/2}$. Moreover,
\bnum
\item $[\Delta_n(\beta)] = \frac{1-q_\beta^{2n}}{1-q_\beta^{2}} [L(\beta)]$,
\item $\Delta_n(\beta)$ is a cyclic module with simple head isomorphic to $L(\beta)$ and socle isomorphic to $q_\beta^{2(n-1)}L(\beta)$,
\item for $n\ge1$, we have
\eqn
&&\Ext_{R(\beta)}^k( \Delta_n(\beta), L(\beta) )
\simeq \left\{
 \begin{array}{ll}
    q_\beta^{-2n} \bR & \hbox{ if } k =1, \\
  0 & \hbox{ if } k \ge 2.
  \end{array}
 \right.
\eneqn
\end{enumerate}
\end{lemma}

Define the {\em root module}
\begin{align*}
\Delta(\beta) := \prolim[n] \Delta_n(\beta).
\end{align*}

\Th[{\cite[Th.\ 3.3]{BKM14}}] \label{Thm: Root module}
There is a short exact sequence
$$
0 \To q_\beta^{2} \Delta(\beta)\To[z_\beta]\Delta(\beta)\To L(\beta) \To 0.
$$
Moreover,
\bnum
\item $\Delta(\beta)$ is a cyclic module with $[\Delta(\beta)] = [L(\beta)]/(1-q_\beta^2)$,
\item $L(\beta)$ is the head of $\Delta(\beta)$,

\item $\tEnd_{R(\beta)}(\Delta(\beta)) \simeq \bR[{z_\beta}]$.
\end{enumerate}
\enth

\Cor [{\cite[Cor.\ 3.5]{BKM14}}] \label{Cor: Root module}
Any finitely generated graded $R(\beta)$-module with all simple subquotients isomorphic to $L(\beta)$ \ro up to a grading shift\/\rf\ is a finite direct sum of grade-shifted copies
of the indecomposable modules $\Delta_n(\beta)$ $(n \ge 1)$ and $\Delta(\beta)$.
\encor

\Prop\label{prop: Root module}
For any $\beta\in\pr$,
$(\Delta(\beta),z_\beta)$ is a strong affinization of $L(\beta)$.
\enprop
\Proof
We can easily check that
conditions~\eqref{hyp:1} and \eqref{hyp:2} in Definition~\ref{def:aff}
are satisfied.

We shall show \eqref{hyp:nd}
by induction on $\Ht{\beta}$.
If $\beta$ is a simple root, then \eqref{hyp:nd} is obvious.
Assume that $\height{\beta}>1$. Then, by \cite[Lemma 4.9, Theorem 4.10]{BKM14},
there exist $\al$, $\gamma\in\pr$ such that
$\al+\gamma=\beta$ and there exists an exact sequence
$$0\to q^{-(\al,\gamma)}\Delta(\gamma)\conv\Delta(\al)
\To[\vphi] \Delta(\al)\conv\Delta(\gamma)\To{[1+p]}\Delta(\beta)\To0.$$
Here $p$ is some non-negative integer and $[1+p]$ is the $q$-integer
with respect to the short root.
Moreover $\vphi$ is given by
\eq
&&\vphi(u\tens v)=\tau_{w[m,n]}(v\tens u)\label{eq:vphiw}
\eneq
for any $u\in\Delta(\gamma)$ and $v\in\Delta(\al)$.
Here $m=\Ht{\al}$ and $n=\Ht{\gamma}$.

By the induction hypothesis, $(\Delta(\al),z_\al)$ and $(\Delta(\gamma),z_\gamma)$
are affinizations.
By \eqref{eq:vphiw},
we see that $\vphi$ commutes with $z_\al$ and $z_\gamma$.
Then $\vphi=a(z_\al,z_\gamma)\Rm_{\Delta(\gamma),\Delta(\al)}$ for some
$a(z_\al,z_\gamma)\in\cor[z_\al,z_\gamma]$ by Proposition~\ref{prop:MNNM}.

Note that $\ap\vert_{\Delta(\al)\circ \Delta(\gamma)}=
\bl \ap\vert_{\Delta(\al)}\br\conv\bl\ap\vert_{\Delta(\gamma)}\br$,
and $\ap\vert_{\Delta(\al)}=c_1z_{\al}^{s_1}$ and $\ap\vert_{\Delta(\gamma)}=c_2z_{\gamma}^{s_2}$ for $c_1,c_2\in\cor^\times$ and $s_1,s_2\in\Z_{\ge0}$.
Hence, if \eqref{hyp:nd} failed, then
$(z_\al z_\gamma)^s\vert_{\Delta(\beta)}=0$ for some $s>0$.
Then we have
$$(z_\al z_\gamma)^s \Delta(\al)\conv\Delta(\gamma)
\subset \Im(\vphi)\subset \Im\bl\Rm_{\Delta(\gamma),\Delta(\al)}\br.$$
Let us take $f(z_\al,z_\gamma)\in\cor[z_\al,z_\gamma]$
such that
$\Rm_{\Delta(\gamma),\Delta(\al)}\Rm_{\Delta(\al),\Delta(\gamma)}
=f(z_\al,z_\gamma)\id_{\Delta(\al)\conv\Delta(\gamma)}$.
Then we have
$$(z_\al z_\beta)^s \Im\bl\Rm_{\Delta(\al),\Delta(\gamma)}\br
\subset f(z_\al,z_\gamma)\Delta(\gamma)\conv\Delta(\al).$$
By Lemma~\ref{lem:monic}, we have
$ f(z_\al,0)\not=0$ and $ f(0,z_\gamma)\not=0$,
which implies
$$ \Im\bl\Rm_{\Delta(\al),\Delta(\gamma)}\br
\subset f(z_\al,z_\gamma)\Delta(\gamma)\conv\Delta(\al).$$
Therefore $f(z_\al,z_\gamma)^{-1}\Rm_{\Delta(\al),\Delta(\gamma)}$ is well defined, which implies that $f$ is an invertible element of $\cor$.
Hence $\Rm_{\Delta(\al),\Delta(\gamma)}$ is an isomorphism.
Then $L(\al)\conv L(\gamma)$ is simple, which is a contradiction.
\QED

Note that
$[L(\beta)]\in[R\gmod]\simeq\Aq$
coincides with the dual PBW vector $E^*(\beta)\in\Aq$.
It is known that $\{E^*(m_{1},\ldots, m_{N})\}_{(m_{1},\ldots, m_{N})\in\Z_{\ge0}^N}$ is a basis of $\Aq$, which is called the
{\em dual PBW basis}.
Here, we set
$$E^*(m_{1},\ldots, m_{N})\seteq
\Bigl( q_{\beta_1} ^{m_1(m_1-1)/2}E^*(\beta_1)^{m_1}\Bigr)\cdots
\Bigl( q_{\beta_N} ^{m_N(m_N-1)/2}E^*(\beta_N)^{m_N}\Bigr)$$
with $\beta_{N-k+1}\seteq r_{i_1}\cdots r_{i_{k-1}}\al_{i_k}$ 
and
$q_{\beta} = q^{(\beta, \beta)/2}$ ($k=1,\ldots,N$).
On the other hand, $E(\beta)=\dfrac{E^*(\beta)}{(E^*(\beta),E^*(\beta))}$
is called the PBW vector and
$$\{E(\beta_1)^{(m_1)}\cdots E(\beta_N)^{(m_N)}\}_{(m_1,\ldots, m_N)\in\Z_{\ge0}^N}$$
is the basis of $\Uam$ called the {\em PBW basis}.
Here $E(\beta)^{(m)}=\dfrac{E(\beta)^m}{[m]_i!}$ with
$i\in I$ such that $(\beta,\beta)=(\al_i,\al_i)$.
Note that the PBW basis and the dual PBW basis are dual to each other.

\section{The duality functor}

\subsection{Duality data}
Let $R$ be the quiver Hecke algebra associated with a generalized Cartan matrix $\cmA$ and polynomials $\qQ_{i,j}(u,v)$.

\Def \label{Def: duality}
Let $J$ be a finite index set.
We say that $\mathcal{D}=\{\beta_j, M_j, \zj,\rr,\Rmat_{j,k}\}_{j,k\in J}$
is a {\em duality datum} if $\beta_j\in\rtlp\setminus\{0\}$,
$M_j \in \Modg(R(\beta_j))$ and
homogeneous homomorphisms
\begin{equation}\label{Eq: endmorphims}
\begin{aligned}
&\mathsf{z}_j \in \tEnd_{R(\beta_j)}(M_j), \quad \mathsf{r}_j \in \tEnd_{R(2\beta_j)}(M_j \circ  M_j), \\
&\mathsf{R}_{j,k} \in \tHom_{R(\beta_j + \beta_k )}(M_j \circ  M_k, M_k \circ  M_j)\quad
\text{for $j,k\in J$}
\end{aligned}
\end{equation}
satisfy the following conditions:
\begin{enumerate}
\item[$\Fa$] For $j\in J$, $\deg \mathsf{z}_j \in 2\Z_{>0}$.
In addition, $M_j$ is a finitely generated free module
over the polynomial ring $\bR[\mathsf{z}_j]$.
\item[$\Fb$] For $j\in J$, we have $\mathsf{r}_j \in \tEnd_{R(2\beta_j)}(M_j \circ  M_j)_{- \deg \mathsf{z}_j}$ and
$$ \mathsf{R}_{j,j} = \bl\mathsf{z}_j \conv {M_j}-{M_j}  \conv \mathsf{z}_j\br
\, \mathsf{r}_j
+ \id_{M_j\circ M_j}. $$
\item[$\Fc$] For $k,l \in J$,
\bna
\item $ (\mathsf{z}_l \conv M_k)\, \mathsf{R}_{k,l}  =
\mathsf{R}_{k,l}\,( M_k \circ \mathsf{z}_l) $ in
$ \tHom_{R(\beta_k + \beta_l)} (M_k \conv M_l, M_l \conv M_k  )$,
\item $(M_l \conv \mathsf{z}_k)\, \mathsf{R}_{k,l}  = \mathsf{R}_{k,l}\,
 (\mathsf{z}_k \conv M_l) $ in
$\tHom_{R(\beta_k + \beta_l)} (M_k \conv M_l, M_l \conv M_k  )$.
\end{enumerate}
\item[$\Fd$] There exist polynomials $\qQ^{\mathcal{D}}_{k,l}(u,v) \in \bR[u,v]$ $(k,l\in J )$
such that
\bna
\item $\qQ^{\mathcal{D}}_{k,k}(u,v)=0$, and
$\qQ^{\mathcal{D}}_{k,l}(u,v)$ $(k\ne l)$ is of the form
$$ \sum_{ \deg \mathsf{R}_{k,l} + \deg \mathsf{R}_{l,k} - p \deg \mathsf{z}_k  - q \deg \mathsf{z}_l=0  } t_{k,l;p,q} u^p v^q, $$
where $t_{k,l;\; (\deg \mathsf{R}_{k,l} + \deg \mathsf{R}_{l,k})/\deg \mathsf{z}_k ,\, 0} \in \bR^\times$,
\item $\qQ^{\mathcal{D}}_{k,l}(u,v) = \qQ^{\mathcal{D}}_{l,k}(v,u) $,
\item $ \mathsf{R}_{l,k} \mathsf{R}_{k,l} = \left\{
                                                     \begin{array}{ll}
                                                       1 & \hbox{ if } k=l, \\
                                                       \qQ^{\mathcal{D}}_{k,l}(\mathsf{z}_k \conv M_l , M_k\conv \mathsf{z}_l ) & \hbox{ if } k \ne l.
                                                     \end{array}
                                                   \right.
$
\end{enumerate}
\item[$\Fe$]
For any $j,k,l \in J$,
$$
(\mathsf{R}_{k,l} \conv M_j) ( M_k\conv\mathsf{R}_{j,l})
(\mathsf{R}_{j,k} \conv M_l) =
( M_l\conv \mathsf{R}_{j,k})(\mathsf{R}_{j,l} \conv M_k) ( M_j\conv \mathsf{R}_{k,l})
$$
holds in $\tHom_{R(\beta_j + \beta_k + \beta_l)}
( M_j \conv M_k \conv M_l, M_l \conv M_k \conv M_j )$.
\end{enumerate}
\edf

For simplicity, we write shortly $\{ M_j, \zj,\Rmat_{j,k}\}_{j,k\in J} $ for $\{\beta_j, M_j, \zj,\rr,\Rmat_{j,k}\}_{j,k\in J}$ if there is no afraid of confusion.

\medskip
We now construct a Cartan datum corresponding to the duality datum
$\mathcal{D}$ as follows.
Let $\{\alD_j\}_{j\in J}$ be the simple roots.
Then we define the weight lattice $\PoD$ by $\PoD=\rtlD\seteq\soplus_{j\in J}\Z\alD_j$,
and define a symmetric bilinear form on $\PoD$ by
\eq
&&
(\alD_j,\alD_k) =
\bc
\deg \zj&\text{if $j=k$,}\\[1ex]
-\bl\deg\qQD_{j,k}(\zj,\zj[k])\br/2=-\bl\deg \Rmat_{j,k}+\deg \Rmat_{k,j}\br/2
&\text{otherwise.}
\ec
\eneq
Define $\hD_j$ by \eqref{item:inner}\;\eqref{item:h} in Definition~\ref{def:cartan}.
Then the corresponding generalized Cartan matrix
$ \cmA^{\mathcal{D}}\seteq( a^{\mathcal{D}}_{jk})_{j,k\in J}$ is given by
$a^{\mathcal{D}}_{jk}=\dfrac{2(\alD_j,\alD_k)}{(\alD_j,\alD_j)}$.
Since $\qQD_{j,k}(\zj,0)\in \cor^\times
\zj^{-a^{\mathcal{D}}_{jk}}$ for $j\not=k$,  $-a^{\mathcal{D}}_{jk}$
 is a non-negative integer. Therefore, $\cmA^{\mathcal{D}}$ is a
generalized Cartan matrix.
We now define $R^{\mathcal{D}}$
as the quiver Hecke algebra   corresponding to the datum
$\{\qQD_{j,k}\}_{j,k\in J}$.

We now have two different quiver Hecke algebras $R$ and $R^{\mathcal{D}}$. To distinguish them,
we write
$$\text{$\xd_k$ ($1 \le k \le \Ht{\gamma})$ and  $\td_l $
($1 \le l \le \Ht{\gamma}-1$)} $$
for the generators $x_k$ ($1 \le k \le \Ht{\gamma}$)
and $\tau_l$ ($1 \le j \le \Ht{\gamma}-1$) of $R^{\mathcal{D}}(\gamma)$
($\gamma\in \rtlpD$).

The $\Z$-grading on $R^{\mathcal{D}}(\gamma) $ is given as follows:
\begin{align*}
\deg(e(\mu))=0, \quad \deg(e(\mu) \xd_k )= \deg \zj[{\mu_k}] , \quad
\deg(e(\mu) \td_l )=
\bc
 -\deg \zj[\mu_l]& \text{if $\mu_{l} = \mu_{l+1}$,} \\
\deg \Rmat_{\mu_{l}, \mu_{l+1}} & \text{if $\mu_{l} \ne \mu_{l+1}$,}
\ec
\end{align*}
which is well-defined (see Definition \ref{dfn: grading}).

Let $\gamma \in {\rlQ}_+^{\mathcal{D}}$ with $m = \Ht{\gamma} $, and define
$$ \Delta^{\mathcal{D}}( \gamma )\seteq\soplus_{\mu \in J^{ \gamma}}
\Delta^{\mathcal{D}}_{\;\mu}, $$
where
$$ \Delta^{\mathcal{D}}_{\;\mu} \seteq M_{\mu_1} \conv M_{\mu_2} \conv \cdots \conv
M_{\mu_m} \quad \text{for $\mu = (\mu_1, \mu_2, \ldots, \mu_m) \in J^{\gamma}$.} $$
Let
$$\phi\cl \rtlD\to\rtl$$
be the linear map defined by $\phi(\alD_j)=\beta_j$ for $j\in J$.
Then, it is clear that $\Delta^{\mathcal{D}}( \gamma)$ is a left
$R(\phi(\gamma))$-module.

We define a right $R^{\mathcal{D}}(\gamma)$-module structure on $\Delta^{\mathcal{D}}(\gamma)$ as follows:
\bna
\item
$e(\mu)$ is the projection to the component $\Delta^{\mathcal{D}}_{\;\mu}$,
\item
the action of $\xd_k$ on $\Delta^{\mathcal{D}}_{\;\mu}$
is given by $M_{\mu_1}\conv \cdots \conv M_{\mu_{k-1}}\conv
\zj[{\mu_k}]\conv M_{\mu_{k+1}}\conv\cdots\conv M_{\mu_{m}}$,
\item if $\mu_{k}\not=\mu_{k+1}$,
the action of $\td_k$ on $\Delta^{\mathcal{D}}_{\;\mu}$ is given by
$$M_{\mu_1}\conv \cdots \conv M_{\mu_{k-1}}\conv
\Rmat_{\mu_k,\,\mu_{k+1}}\conv M_{\mu_{k+2}}\conv\cdots \conv M_{\mu_{m}},$$
\item
if $\mu_{k}=\mu_{k+1}$, the action of $\td_k$ on $\Delta^{\mathcal{D}}_{\;\mu}$
is given by
$$M_{\mu_1}\conv \cdots \conv M_{\mu_{k-1}}\conv
\rr[{\mu_k}]
\conv M_{\mu_{k+2}}\conv\cdots \conv M_{\mu_{m}}.$$
\ee

\Th \label{Thm: bimodule}
The right $R^{\mathcal{D}}(\gamma)$-module structure on $\Delta^{\mathcal{D}}( \gamma )$
is well-defined.
\enth
\Proof
Since the proof is easy and similar to the arguments in \cite{K^3}, we omit it.
\QED

By the construction, the right $R^{\mathcal{D}}(\gamma)$-module action commutes with the left $R(\phi(\gamma))$-module action, which means that
\begin{align*}
\text{ $\Delta^{\mathcal{D}}(\gamma)$ has an $( R(\phi(\gamma)), R^{\mathcal{D}}(\gamma) )$-bimodule structure.}
\end{align*}
We now define the functor
\begin{align*}
\F^{\mathcal{D}}_{\gamma} \cl \Modg(R^{\mathcal{D}}(\gamma))\To
\Modg\bl R(\phi(\gamma))\br
\end{align*}
by
$$ \F^{\mathcal{D}}_{\gamma}(M):= \Delta^{\mathcal{D}}(\gamma) \otimes_{R^{\mathcal{D}}(\gamma)} M. $$
Set
\begin{align*}
\F^{\mathcal{D}}= \bigoplus_{\gamma \in \rtlpD} \FD_{\gamma}.
\end{align*}

For $j\in J$, we shall write
$\LD(j)$ for the simple $R^\D(\alD_j)$-module $R^\D(\alD_j)/ R^\D(\alD_j)\xd_1$.

\Th\label{Thm: main thm} Let  $\mathcal{D}=\{\beta_j, M_j, \zj,\rr,\Rmat_{j,k}\}_{j,k\in J}$ be a duality datum.
\bnum
\item The functor
$$ \F^{\mathcal{D}} : \Modg(R^{\mathcal{D}})\To\Modg(R) $$
is a tensor functor.
\item
For $j\in J$,
$$\text{
$\F^{\mathcal{D}} \bl R^{\mathcal{D}}(\alD_j)\br\simeq \Ma_j$ and
$\F^{\mathcal{D}} (\LD(j)) \simeq \Ma_j / \mathsf{z}_j \Ma_j$}. $$

\item If $\cmA^\mathcal{D}$ is of finite type, then the functor $\F^{\mathcal{D}}$ is exact.
\item If a graded $R^{\mathcal{D}}(\gamma)$-module $L$ is finite-dimensional, then so is $\F^{\mathcal{D}}(L)$. Thus, we have the induced functor
$$ \F^{\mathcal{D}} \cl R^\mathcal{D}\gmod \longrightarrow
R\gmod.$$
\end{enumerate}
\enth
\begin{proof}
Since the proof is easy and similar to \cite{K^3}, we omit it.
\end{proof}

\subsection{Construction of duality data from affinizations}

Let $J$ be a finite index set.
Let $ \{\beta_j, \Ma_j, \zj\}_{j\in J}$ be a datum such that
\bna
\item $\beta_j\in\rtlp\setminus\{0\}$,
\item $(\Ma_j,\zj)$ is an even affinization of a real simple $R(\beta_j)$-module
$\bM_j\seteq \Ma_j/\zj \Ma_j$.
\setcounter{mycounter}{\value{enumi}}
\ee
Then we take $\Rmat_{j,k}$ as follows:
\bna\setcounter{enumi}{\value{mycounter}}
\item $\Rmat_{j,k}=\Rm_{\Ma_j,\Ma_k}$.
Furthermore, we normalize $\Rmat_{j,j}$
such that $\Rmat_{j,k}\vert_{\zj=\zj[k]=0}=\id_{\bM_j\circ\bM_j}$ when $j=k$.
\ee

Then, Proposition~\ref{prop:MM}
implies that
\eq
\rr\seteq(\zj\conv \Ma_j-\Ma_j\conv \zj)^{-1}\bl \Rmat_{j,j}-\id_{\Ma_j\circ \Ma_j}\br\label{def:rr}
\eneq
is a well-defined endomorphism of $\Ma_j\conv\Ma_j$.

Note that for any $ \{\beta_j,  \Ma_j, \zj\}_{j\in J}$ satisfying
(a) and (b), we can always choose  $\Rmat_{j,k}$'s.
Moreover, $\Rmat_{j,j}$ is unique and
$\Rmat_{j,k}$ ($j\not=k$) is unique up to constant multiple.

\Th \label{Thm: duality and affinization}
Under the above assumptions {\rm (a), (b), (c)}, we have the following.
\bnum
\item The datum $\mathcal{D}=\{\beta_j, \Ma_j, \zj,\rr,\Rmat_{j,k}\}_{j,k\in J}$
is a duality datum.
\item \label{item:finite} Assume that $\cmA^\mathcal{D}$ is of finite type.
Then, we have the following.
\bna
\item
$\F^\mathcal{D}(M)$ is either a simple module or vanishes
for any simple  $R^\mathcal{D}$-module $M$.

Moreover, if $M$ is a real simple module and $\F^{\mathcal{D}}(M)$ is non-zero, then $\F^{\mathcal{D}}(M)$ is real.

\item \label{item:aff}
Let $(\Na, \zN)$ is an affinization of a simple $R^{\mathcal{D}}$-module $\bN$.
If $\F^{\mathcal{D}}(\bN)$ is simple, then $(\F^{\mathcal{D}} (\Na), \F^{\mathcal{D}}(\zN))$ is an affinization of $\F^{\mathcal{D}}(\bN)$.

\item  Let $M$ and $N$ be simple $R^\mathcal{D}$-modules, and assume that
one of them is real and also admits an affinization.
Then $\F^\mathcal{D}(M\hconv N)$ is zero or isomorphic to
$\F^\mathcal{D}(M)\hconv \F^\mathcal{D}(N)$.
\ee
\ee
\enth

\Proof (i)\
Let us prove that $\mathcal{D}$ is a duality datum.
Since axioms $\Fa$--$\Fd$ are obvious,
we only gives the proof of the braid relation $\Fe$:
\eq
&&\Rmat_{jk}\circ\Rmat_{ik}\circ\Rmat_{ij}
=\Rmat_{ij}\circ\Rmat_{ik}\circ\Rmat_{jk}\label{eq:braidR}
\eneq
as a morphism from
$\Ma_i\conv \Ma_j\conv \Ma_k\to \Ma_k\conv \Ma_j\conv \Ma_i$ for $i,j,k\in J$.
By the definition, we have
$R_{\Ma_i,\Ma_j}=a(\zj[i],\zj)\Rmat_{i,j}$ for a non-zero polynomial $a(\zj[i],\zj)$.
The R-matrices $R_{\Ma_i,\Ma_j}$ satisfy the braid relation:
$$R_{\Ma_j,\Ma_k} \circ R_{\Ma_i,\Ma_k}\circ R_{\Ma_i,\Ma_j}=
R_{\Ma_i,\Ma_j}\circ R_{\Ma_i,\Ma_k}\circ R_{\Ma_j,\Ma_k}.$$
The calculation
\begin{align*}
R_{\Ma_j,\Ma_k} \circ R_{\Ma_i,\Ma_k}\circ R_{\Ma_i,\Ma_j}
&=a(\zj,\zj[k])\Rmat_{j,k}\circ a(\zj[i],\zj[k])\Rmat_{i,k}
\circ a(\zj[i],\zj)\Rmat_{i,j}\\
&=a(\zj,\zj[k])a(\zj[i],\zj[k])a(\zj[i],\zj)\Rmat_{j,k}\circ\Rmat_{i,k}
\circ\Rmat_{i,j},
\end{align*}
and a similar calculation for $R_{\Ma_i,\Ma_j}\circ R_{\Ma_i,\Ma_k}\circ R_{\Ma_j,\Ma_k}$
show that
$$a(\zj,\zj[k])a(\zj[i],\zj[k])a(\zj[i],\zj)
\Bigl(\Rmat_{jk}\circ\Rmat_{ik}\circ\Rmat_{ij}
-\Rmat_{ij}\circ\Rmat_{ik}\circ\Rmat_{jk}\Bigr)=0.$$
Hence we obtain \eqref{eq:braidR}.

\medskip\noi
(ii) (a)\hs{1.5ex}Let us prove that
$\F^\mathcal{D}(M)$ is a simple module or zero
for a simple $R^\mathcal{D}(\gamma)$-module $M$ by induction on
$\height{\gamma}$.
Let us assume $M\simeq N\hconv \LD(j)$ for some $j\in J$ and
a simple $R^\mathcal{D}(\gamma-\alD_j)$-module $N$.
By the induction hypothesis, $\F^\mathcal{D}(N)$ is a simple module or zero.
Let $r\col N\conv \LD(j)\to \LD(j)\conv N$ be a non-zero homomorphism
of $R^\mathcal{D}(\gamma)$-modules.
Then $\Im(r)$ is isomorphic to $N\hconv \LD(j)$.
Since $\F^\mathcal{D}$ is exact, $\F^\mathcal{D}(\Im(r))\simeq
\Im\bl \F^\mathcal{D}(r)\br\simeq \F^\mathcal{D}(M)$.
If $\F^\mathcal{D}(N)\simeq 0$, then
$\F^\mathcal{D}(M)\simeq 0$.
Assume that $\F^\mathcal{D}(N)$ is a simple module.
Then $\Im\bl \F^\mathcal{D}(r)\br$ is isomorphic to
$\F^\mathcal{D}(N)\hconv
\F^\mathcal{D} (\LD(j))$ or $0$ according that
$\F^\mathcal{D}(r)$ is non-zero or zero by Proposition~\ref{Prop: simple-unique}.

If $M$ is real simple and $\F^\mathcal{D}{M}$ is simple, then
$(\F^\mathcal{D}{M})\conv (\F^\mathcal{D}M)\simeq \F^\mathcal{D}(M\conv M)$
is simple and hence $\F^\mathcal{D}{M}$ is real.

Thus we obtain  (ii) (a).

\medskip\noi
(ii) (b)\
Let $\bN$ be a simple $R^{\mathcal{D}} (\gamma)$-module and set $m = \Ht{\gamma}$.
We put $ \Na_\F = \F^{\mathcal{D}}(\Na)$ and $z_\F = \F^{\mathcal{D}}(\zN)$.
Applying the functor $\F^{\mathcal{D}}$ to the exact sequence
$$ 0 \longrightarrow \Na \To[\;\zN\; ] \Na \longrightarrow \bN \longrightarrow 0,  $$
we obtain the exact sequence
$$ 0 \longrightarrow \Na_\F \To[\;z_\F\;] \Na_\F \longrightarrow \F^{\mathcal{D}}(\bN) \longrightarrow 0.  $$
Thus, we have an injective homogeneous endomorphism $z_\F$ of $\Na_\F$ and $ \Na_\F / z_\F \Na_\F \simeq  \F^{\mathcal{D}}(\bN) $.
Since $\Ma_j$ is a finitely generated $R(\beta_j)$-module for any $j$ by
Lemma~\ref{lem:fg}, $\Na_\F$ is a finitely generated graded $R$-module
and $\F^{\mathcal{D}}(\bN)$ is a finite-dimensional $R$-module.
Hence, condition \eqref{hyp:1} of Definition \ref{def:aff} holds
(see Remark~\ref{rem:aff} \eqref{it: affeq}).

Let us show \eqref{hyp:nd} of Definition \ref{def:aff}.
Let $i\in I$.
By  Lemma \ref{lem:affweek} \eqref{cor:MNNM},
for any $j\in J$, there exist
$d_j\in\Z_{\ge0}$ and  $c_j\in\cor^\times$
such that
$ \ap |_{\Ma_{j}}=c_j \zj^{d_j}$.
Since
$ \ap |_{\Ma_{\mu_1} \circ \cdots \circ \Ma_{\mu_m}  } =  (\ap|_{\Ma_{\mu_1}}) \conv \cdots \conv
(\ap|_{\Ma_{\mu_m}}  )=\prod_{k=1}^mc_{\mu_k}(\xd_k)^{d_{\mu_k}} $, we obtain
\begin{align*}
\ap\vert_{ \Na_\F} &= \sum_{\mu \in J^\gamma} \Bigl(\ap \vert_{
\Ma_{\mu_1} \circ \cdots \circ \Ma_{\mu_m} }\Bigr) \tens_{R^{\mathcal{D}}(\gamma)} \Na   \\
&= \sum_{\mu \in J^\gamma}  ( \Ma_{\mu_1} \conv \cdots \conv \Ma_{\mu_m} )
\tens_{R^{\mathcal{D}}(\gamma)} \Bigl( e(\mu) c (\xd_1)^{d_{\mu_1}}
\cdots (\xd_m)^{d_{\mu_m}} \Bigr)\big\vert _\Na   \\
&= \sum_{\mu \in J^\gamma}  ( \Ma_{\mu_1} \conv \cdots \conv \Ma_{\mu_m} )
\tens_{R^{\mathcal{D}}(\gamma)} \Bigl(c e(\mu)
\prod_{j\in J}\bl\smash{\prod_{\substack{k\in[1,m],\;\mu_k=j}}}(\xd_k)^{d_j} \br
\Bigr)\big\vert _\Na   \\
&=   \Delta^{\mathcal{D}}(\gamma) \tens_{R^{\mathcal{D}}(\gamma)} (  c \prod_{j\in J} \ap[j]^{d_{j}} ) \big\vert_ \Na
\end{align*}
with $c=\prod_{k=1}^mc_{\mu_k}$ which does not depend on $\mu\in J^\gamma$.
 Therefore, condition \eqref{hyp:nd}
of Definition \ref{def:aff} holds.

\medskip\noi
(ii) (c) immediately follows from (a) and
the epimorphism
$$ \F^\mathcal{D}(M)\conv \F^\mathcal{D}(N)
\epito \F^\mathcal{D}(M\hconv N)$$
because $M\hconv N$ is simple.
\QED

\section{Examples} \label{Sec: Ex}

Let $\g$ be a Kac-Moody Lie algebra associated with a Cartan matrix $\cmA$
of finite type. Suppose that
\begin{equation} \label{Eq: condiotions on beta}
\left\{\parbox{72ex}{
\bna
\item
$\{ \beta_j\}_{ j \in J }$  is a family of elements of $\pr$,
which is linearly independent in $\rlQ$,
\item
$\beta_j - \beta_k \notin \Phi$ for any $j,k \in J$, where $\Phi$ is the set of roots of $\g$.
\ee
}\right.
\end{equation}
Let $\overline{\g}$ be the Lie subalgebra of $\g$ generated by
the root vectors of weight $\beta_j$ and $-\beta_j$ (cf.\ \cite[Th.\ 1.1]{Naito91}).
Then $\overline{\g}$ is a Kac-Moody Lie algebra associated to
\begin{align} \label{Eq: bar A}
\overline{\cmA}:= (\overline{a}_{j,k})_{j,k\in J}\quad
\text{with $\overline{a}_{j,k}\seteq 2(\beta_j, \beta_k)/(\beta_j,\beta_j)$.}
\end{align}
We have an injective algebra homomorphism
\begin{align} \label{Eq: embedding}
\xymatrix@C=2.5ex{U^-(\overline{\g})\;\; \ar@{>->}[r]& U^-(\g).}
\end{align}

Choosing a convex order of the set $\pr$
of positive roots,
let $\bl\Delta(\beta_j),\zj\br$ be the affinization of $L(\beta_j)$
given in Proposition~\ref{prop: Root module}.
Then, we have the duality datum
$$\mathcal{D} := \{ \Delta(\beta_j), \mathsf{z}_j, \mathsf{R}_{k,l}  \}_{j,k,l\in J  } .$$
Let $\g^{\mathcal{D}}$ be the Kac-Moody Lie algebra associated with
$\cmA^{\mathcal{D}}$.
Suppose that $\cmA^{\mathcal{D}}$ is of finite type.
Then, the functor $\F^{\mathcal{D}}$ is exact, and
gives a $\Z[q^{\pm1}]$-algebra homomorphism:
$$[R^{\mathcal{D}}\gmod]\To {[R\gmod]}$$
which gives the $\Z[q^{\pm1/2}]$-algebra homomorphism (see Corollary~\ref{cor:twist}):
\eq
&&\Aq[{(\g^{\mathcal{D}})^+}]_\fc\to \Z[q^{\pm1/2}]\otimes_{\Z[q^{\pm1}]}\Aq.
\label{eq:Aqhom}
\eneq
sending $f_j$ to the dual PBW generator $E^*(\beta_j)$
corresponding to $ [\Delta(\beta_j)]$.
Here $\fc$ is the bilinear form on $\rtl^{\D}$
given by $\fc(\al^{\D}_j,\al^{\D}_k)=\frac{1}{2}(\deg \mathsf{R}_{k,j}-\deg \mathsf{R}_{j,k})$.

By applying the exact functor $ \Q(q^{1/2})\tens_{\Q[q^{\pm1/2}]}\scbul$ to \eqref{eq:Aqhom},
we obtain a $\Q(q^{1/2})$-algebra homomorphism
\eq\Uqm[{\g^\D}]_\fc\To \Q(q^{1/2})\tens_{\Q(q)}\Uqm.
\label{mor:mD}
\eneq
Set $c_\beta\seteq(E^*(\beta),E^*(\beta))^{-1}$.
Then $E(\beta)=c_{\beta}E^*(\beta)$ is the PBW vector corresponding to $\beta\in\pr$.
Let $\psi$ be the algebra automorphism of
$\Uqm[{\g^\D}]_\fc$ sending $f_j$ to $c_{\beta_j}f_j$.
Then the composition
$$\xymatrix{
\Uqm[{\g^\D}]_\fc\ar[r]_-\psi^\sim&\Uqm[{\g^\D}]_\fc\ar[r]&
\Q(q^{1/2})\tens_{\Q(q)}\Uqm}$$
sends $f_j$ to $E(\beta_j)$.
Since $\deg \z=(\beta_j,\beta_j)$ by Theorem~\ref{Thm: Root module},
the above homomorphism sends  the divided power $f_j^{(m)}$ to the divided power
$E(\beta_j)^{(m)}$.
Moreover, the $f_j^{(m)}$'s generate the $\A$-algebra
$\Uam[{\g^\D}]_\fc$,
and the $E(\beta_j)^{(m)}$'s are contained in $\Uam$.
Hence we obtain the algebra homomorphism
\eq\Uam[{\g^\D}]_\fc\To \Q[q^{\pm1/2}]\tens_{\Q[q^{\pm1}]}\Uam.
\label{mor:UqmD}
\eneq
Taking the classical limit $q^{1/2}=1$, we obtain the induced algebra homomorphism
\eq U^- ( \g^{\mathcal{D}}  ) \longrightarrow U^- ( \g  )
\label{mor:Fdlass}
\eneq
sending $f_j$ to the root vector corresponding to $-\beta_j$ for $j\in J$.

\Prop\label{prop:injcr} If
$\cmA^{\mathcal{D}} = \overline{\cmA}$,
then the morphism induced by $\FD$
$$[R^{\mathcal{D}}\gmod]\To {[R\gmod]}$$
is injective.
In particular $\FD$ sends the simple $\RD$-modules to simple $R$-modules.
\enprop
In such a case,
the functor $\F^{\mathcal{D}}$ categorifies the homomorphism $\eqref{Eq: embedding}$.
\Proof
By the condition, we have
$U^- ( \g^{\mathcal{D}}  ) \simeq U^- (\ol{\g})$.
Hence the map \eqref{mor:Fdlass} is injective, which implies that
\eqref{mor:UqmD} is injective.
Hence
\eqref{mor:mD}  as well as  \eqref{eq:Aqhom} is injective.
\QED

Let us give several examples of such duality data.
\Ex \label{Ex: D}
Let $I = \{1,2,\ldots ,\ell \}$ and $\cmA$ the Cartan matrix of type $A_\ell$.
Hence $(\al_i,\al_j)=2\delta(i=j)-\delta(|i-j|=1)$ for $i,j\in I$.
Let $R$ be the quiver Hecke algebra associated with $\cmA$ and the parameter $\qQ_{i,j}(u,v)$ defined as follows: for $i,j \in I$ with $i < j$,
$$ \qQ_{i,j}(u,v) = \left\{
                           \begin{array}{ll}
                             u-v & \text{ if } j = i+1, \\
                             1 & \text{ otherwise. }
                           \end{array}
                         \right.
$$

Let $J = \{1,2,\ldots ,\ell \}$ and $\beta_1\seteq\alpha_1 + \alpha_2$,
$\beta_j \seteq\alpha_j$ for $j\in J\setminus\{1\}$.
Note that the $\beta_j$'s do not satisfy condition \eqref{Eq: condiotions on beta} (b). We put
\begin{align*}
\Delta(\beta_1) := L(1,2)_{\mathsf{z}_1}, \quad \Delta(\beta_j) := L(j)_{\mathsf{z}_j} \quad (j \in J \setminus \{ 1 \}),
\end{align*}
where $L(1,2) := \bR v$ is
the 1-dimensional $R(\beta_1)$-module with the actions
$$ e(\nu) v = \delta_{\nu, (1,2)}v, \quad x_1 v = x_2 v = \tau_1 v = 0 \qquad \text{ for } \nu \in I^{\alpha_1 + \alpha_2} .  $$
Note that $\deg(\mathsf{z}_j)=2$ for $j\in J$ and the $\Delta(\beta_j)$'s are root modules.
We set $\mathsf{R}_{j,k}= \Rm_{\Delta(\beta_j), \Delta(\beta_k)}$.
By direct computations, the R-matrix $R_{\Delta(\beta_j), \Delta(\beta_k)}$ $(j\ne k)$ is given as follows: for $u \otimes v \in \Delta(\beta_j) \otimes \Delta(\beta_k)$,
\begin{align*}
R_{\Delta(\beta_j), \Delta(\beta_k)} (u \otimes v) =
\left\{
     \begin{array}{ll}
        (\tau_2\tau_1(\mathsf{z}_2 - \mathsf{z}_1) + \tau_1)(v \otimes u)  & \hbox{ if $j=1$ and $k = 2$},  \\
        \tau_2\tau_1 (v \otimes u) & \hbox{ if $j=1$ and $k > 2$}, \\
        \tau_1 \tau_2 (\mathsf{z}_1 - \mathsf{z}_2) (v \otimes u) & \hbox{ if $j=2$ and $k = 1$}, \\
        \tau_1 \tau_2 (v \otimes u) & \hbox{ if $j>2$ and $k =1$}, \\
        \tau_1 (v \otimes u) & \hbox{ otherwise, }
     \end{array}
   \right.
\end{align*}
which yields
$$
\mathsf{R}_{j,k} =
\left\{
  \begin{array}{ll}
    (\mathsf{z}_1 - \mathsf{z}_2)^{-1} R_{\Delta(\beta_j), \Delta(\beta_k)}
& \hbox{ if $ j=2$ and $k=1$},  \\[1ex]
    R_{\Delta(\beta_j), \Delta(\beta_k)} & \hbox{ otherwise,}
  \end{array}
\right.
$$
and
$$ \deg (\mathsf{R}_{j,k} ) = \left\{
                                \begin{array}{ll}
          1 & \hbox{if $ |j-k|=1$ and $(j,k)\not= (2,1)$,}\\
1&\text{$(j,k)=(1,3), (3,1)$,}\\
                                  -1 & \hbox{ if $(j,k) = (2,1)$}, \\
                                  0 & \hbox{ otherwise.}
                                \end{array}
                              \right.
 $$
Thus, we have
$$\cmA^{\mathcal{D}} = \left(
                         \begin{array}{ccccccc}
                           2 & 0 & -1 & 0 &\cdots & 0 &0\\
                           0 & 2 & -1 & 0 & \cdots  &0&0  \\
                           -1 & -1 & 2 & -1 & \cdots &  0&0 \\
                           0 & 0 & -1 & 2 &\cdots & 0 &0 \\
                           \vdots & \vdots & \vdots & \vdots & \ddots & \vdots& \vdots\\
                           0 & 0 & 0 & 0 & \cdots & 2 & -1 \\
                           0 & 0 & 0 & 0 & \cdots & -1 &2
                         \end{array}
                       \right)
$$
which is of type $D_\ell$, i.e., the quiver Hecke algebra $R^{\mathcal{D}}$ is of type $D_\ell$.
Note that $\deg(e(1,2) \td_1 ) = 1$ and $\deg(e(2,1) \td_1 ) = -1$
(see Definition \ref{dfn: grading}).
By Theorem \ref{Thm: duality and affinization}, we have a functor $ \F^{\mathcal{D}} $ between quiver Hecke algebras of type $D_\ell$ and $A_\ell$
such that
$$ \F^{\mathcal{D}}( L^{\mathcal{D}}(j)) \simeq L(\beta_j) \quad  \text{ for $j\in J.$} $$

Let us consider the $R^{\mathcal{D}}$-module
$\LD(1,3) :=\LD(1) \hconv \LD(3)$ and the one-dimensional  $R$-module
$L(1,2,3) := L(1,2) \hconv L(3)$. Applying the functor $\F^{\mathcal{D}}$ to the exact sequence
$$ 0 \rightarrow \LD(1,3) \rightarrow \LD(3)\conv \LD(1) \rightarrow
\LD(1)\conv \LD(3) \rightarrow \LD(1,3) \rightarrow 0, $$
we have
$$ 0 \rightarrow \F^{\mathcal{D}} (\LD(1,3)) \rightarrow L(3)\conv L(1,2) \rightarrow L(1,2)\conv L(3) \rightarrow \F^{\mathcal{D}}(\LD(1,3)) \rightarrow 0. $$
Since $\F^{\mathcal{D}}$ sends a simple module to a simple module or zero by Theorem \ref{Thm: duality and affinization} and $L(3)\conv L(1,2)$
is not isomorphic to $L(1,2)\conv L(3)$, we have
$$ \F^{\mathcal{D}}(\LD(1,3)) \simeq L(1,2,3). $$
Set $\LD(1,3,2)\simeq\LD(1,3)\hconv \LD(2)$, which
is one-dimensional.
It is isomorphic to the image of the composition of
\eqn
&&\LD(1)\conv \LD(3)\conv \LD(2)
\To \LD(1)\conv \LD(2)\conv \LD(3)
\To\LD(2)\conv \LD(1)\conv \LD(3)\\
&&\hs{57ex}\hfill\To\LD(2)\conv\LD(1,3).\eneqn
By applying $\FD$,
we obtain the diagram
\eq
&&\xymatrix{
L(1,2)\conv L(3)\conv L(2)\ar[r]^{f_1}&L(1,2)\conv L(2)\conv L(3)
\ar[r]^{f_2}& L(2)\conv L(1,2)\conv L(3)\ar@{->>}[d]^{f_3}\\
&&L(2)\conv L(1,2,3).}
\label{dia:Rcomp}
\eneq
Hence $\FD(\LD(1,3,2))$ is isomorphic to the image of $f_3f_2f_1$.
Let $u_{1,2}$, $u_2$, and $u_3$ be the generator of
$L(1,2)$, $L(2)$ and $L(3)$, respectively.
Then we have
\eqn
&&f_1(u_{1,2}\tens u_3\tens u_2)=\tau_3(u_{1,2}\tens u_2\tens u_3),\\
&&f_2(u_{1,2}\tens u_2\tens u_3)=\tau_1(u_{2}\tens u_{1,2}\tens u_3).
\eneqn
Therefore, we obtain
\eqn
&&f_2f_1(u_{1,2}\tens u_3\tens u_2)=
\tau_3\tau_1(u_{2}\tens u_{1,2}\tens u_3)
=\tau_1\tau_3(u_{2}\tens u_{1,2}\tens u_3),
\eneqn
which is killed by $f_3$.
Thus $f_3f_2f_1=0$, and hence we conclude
$$\FD(\LD(1,3,2))\simeq 0.$$
Therefore, $\FD$ can send simple modules to zero in this example.
\enEx

\vskip 1em

\section{Further examples for non-symmetric types} \label{Sec: Ex2}

Let $\beta\in\rtlp$ and let
$(\Ma,\z)$ be an affinization of a real simple $R(\beta)$-module $\bM$.
We set $J=\{\zr\}$, $\beta_\zr=\beta$, $\Ma_\zr=\Ma$.
Then
$$\mathcal{D}= \{ \Ma_\zr, \z, \Rm_{\Ma,\Ma}\} $$
is a
duality datum. Then the corresponding simple root $\alD_\zr$ satisfies
$(\alD_\zr,\alD_\zr)=\deg \z$.
Let $\bl \KM(\zr^n),\zj[{\KM(\zr^n)}]\br$ be the affinization of the simple
$\RD(n\alD_\zr)$-module $L(\alD_\zr)^{\circ \,n}$ given in Example~\ref{Ex:Affn}.

Now $\Ma^{\circ\mspace{1mu}n}\seteq\underbrace{\Ma\conv\cdots\conv\Ma}_{\text{$n$ times}}$
has a structure of $(R(n\beta),\RD(n\alD_\zr))$-bimodule.
We set
 $$\Cv_n(\Ma)=\Ma^{\circ\,n}\tens_{\RD(n\alD_\zr)}\KM(\zr^n)
\simeq\FD(\KM(\zr^n)).$$
Then $\zj[{\KM(\zr^n)}]\in\tEnd(\KM(\zr^n))$ induces an endomorphism
$\zj[{\Cv_n(\Ma)}]\in\tEnd(\Cv_n(\Ma))_{n\deg \z}$.
By Theorem~\ref{Thm: duality and affinization}~
\eqref{item:finite}\,\eqref{item:aff},
we obtain the following lemma.
\Lemma $(\Cv_n(\Ma),\zj[{\Cv_n(\Ma)}]) $ is an affinization of
the real simple module $\bM^{\circ\,n}$.
\enlemma

For example,
\begin{align} \label{Eq: def of tilde M}
\Cv_2(\Ma) =
\frac{ \Ma \conv \Ma }{
(\z\conv \Ma + \Ma \conv \z ) (\Ma\conv \Ma) + \rr[](\Ma \conv\Ma)  }  ,
\end{align}
where $\mathsf{r}$ is the endomorphism given in \eqref{def:rr}.
$\zj[{\Cv_2(\Ma)}]$ is the endomorphism induced by $\z\conv \z$.

\medskip
Let $M \in R(\beta)\gmod$ and $N \in R(\beta')\gmod$
be real simple modules.
Suppose that $R(\beta)$ and $R(\beta')$ are symmetric, and
$$
\Rm_{N_{\mathsf{z}'}, M_{\mathsf{t}}} \Rm_{M_{\mathsf{t}}, N_{\mathsf{z}'}} = c (\mathsf{t} - \mathsf{z}')^p
\in \End_{R(\beta+\beta')}( M_{\mathsf{t}} \circ N_{\mathsf{z}'} )
$$
for some $c\in \bR^\times$ and $p \in \Z_{\ge0}$. Set
\begin{align*}
R_1 &:= (\Rm_{M_{\mathsf{t}_1}, N_{\mathsf{z}'}} \conv M_{\mathsf{t}_2})
(M_{\mathsf{t}_1} \conv\Rm_{M_{\mathsf{t}_2}, N_{\mathsf{z}'}}) \in
\tHom_{R(2\beta+\beta')} (M_{\mathsf{t}_1}\conv M_{\mathsf{t}_2} \conv N_{\mathsf{z}'} ,\,
 N_{\mathsf{z}'} \conv M_{\mathsf{t}_1}\conv M_{\mathsf{t}_2}), \\
R_2 &:= (M_{\mathsf{t}_1}\conv\Rm_{N_{\mathsf{z}'}, M_{\mathsf{t}_2}} ) (\Rm_{N_{\mathsf{z}'}, M_{\mathsf{t}_1}} \conv M_{\mathsf{t}_2}) \in \tHom_{R(2\beta+\beta')} ( N_{\mathsf{z}'} \conv M_{\mathsf{t}_1}\conv
M_{\mathsf{t}_2}, M_{\mathsf{t}_1}\conv M_{\mathsf{t}_2} \conv N_{\mathsf{z}'} ).
\end{align*}
Setting $\mathsf{t}_1+\mathsf{t}_2=0$ and
$\mathsf{t}_1\mathsf{t}_2=\widehat{\mathsf{z}}\seteq\zj[\Cv_2(M)]$,
we regard $R_1$ and $R_2$ as homomorphisms in
$\tHom_{R(2\beta+\beta')} \bl \Cv_2(M_\mathsf{z}) \conv N_{\mathsf{z}'}, N_{\mathsf{z}'} \conv
\Cv_2(M_\mathsf{z}) \br$ and
$\tHom_{R(2\beta+\beta')} \bl N_{\mathsf{z}'} \conv \Cv_2(M_\mathsf{z}),\, \Cv_2(M_{\mathsf{z}}) \conv N_{\mathsf{z}'}\br$, respectively.
Then, we have
\begin{equation}  \label{Eq: R1 R2}
\begin{aligned}
R_2 R_1 &= c^2 ( \mathsf{t}_1 - \mathsf{z}')^p ( \mathsf{t}_2 - \mathsf{z}')^p \\
&= c^2 ( \mathsf{t}_1\mathsf{t}_2 - ( \mathsf{t}_1+\mathsf{t}_2 )\mathsf{z}' + \mathsf{z'}^{\,2} )^p \\
&= c^2 ( \widehat{\mathsf{z}} + \mathsf{z'}^{\,2} )^p
\end{aligned}
\end{equation}
in $ \End_{R(2\beta+\beta')}( \Cv_2(M_\mathsf{z})\conv N_{\mathsf{z}'} ) $.

\smallskip
Using $\eqref{Eq: R1 R2}$, one can construct functors $\F^{\mathcal{D}}$ between symmetric and non-symmetric quiver Hecke algebras.
In particular, the functor from type $C_\ell$ (resp.\ $C_\ell^{(1)}$, $A_{2\ell-1}^{(2)}$) to type $A_\ell$ (resp.\ $A_{\ell+1}$, $D_{\ell+1}$) can be constructed.
We give such constructions in the following examples.

\Ex \label{Ex: C}
We take $I$, $\cmA$, and $\qQ_{i,j}(u,v)$ given in Example \ref{Ex: D}.
In particular, $\g$ is of type $A_\ell$.

Let $J = \{1,2,\ldots, \ell \}$ and
$$
\beta_1 = 2\alpha_1,\ \beta_j = \alpha_j\ \text{for $j\in J\setminus\{1\}$}.
$$
Let us denote
\begin{align*}
\Ma_1 = \KM(1^2),\quad \Ma_j = L(j)_{\mathsf{z}_j} \ (j \in J \setminus \{1 \}),
\end{align*}
and $\zj[1] := \zj[K(1^2)]$. Then $\deg \zj[1]=4$ and $\deg \zj=2$ for $j\not=1$.
Note that
$$
\Rm_{L(j)_{\mathsf{z}}, L(k)_{\mathsf{w}} } =
R_{L(j)_{\mathsf{z}}, L(k)_{\mathsf{w}} }.
$$
We put $ \mathsf{R}_{j,k} := R_{\Ma_j, \Ma_k} $.
It follows from $\eqref{Eq: R1 R2}$ that, for $j, k\in J$ with $j < k$,
\begin{align}
\mathsf{R}_{k,j} \mathsf{R}_{j,k} =
\bc
 \mathsf{z}_1  +  \mathsf{z}_2^2  & \hbox{ if } (j,k)=(1,2), \\
 \mathsf{z}_j  -  \mathsf{z}_k
& \hbox{ if $k=j+1$ and $(j,k)\not=(1,2)$,}  \\
                                        1 & \hbox{ otherwise.}
   \ec
\label{eq:AD}
\end{align}

We now set
$$ \mathcal{D} = \{ \Ma_j, \zj, \mathsf{R}_{j,k} \}_{j,k\in J}. $$
Then $\mathcal{D}$ is a dual datum, and \eqref{eq:AD} implies that
$$\cmA^{\mathcal{D}} = \left(
                         \begin{array}{cccccc}
                           2 & -1 & 0  &\cdots & 0 &0\\
                           -2 & 2 & -1  & \cdots  &0&0  \\
                           0 & -1 & 2  & \cdots &  0&0 \\
                           \vdots &  \vdots & \vdots & \ddots & \vdots& \vdots\\
                           0 & 0 &  0 & \cdots & 2 & -1 \\
                           0 & 0 &  0 & \cdots & -1 &2\\
                         \end{array}
                       \right)
$$
which is of type $C_\ell$. Therefore, we have the quiver Hecke algebra $R^{\mathcal{D}}$ of type $C_\ell$ and
the functor $ \F^{\mathcal{D}}  $ from
the category of modules over  quiver Hecke algebras of type $C_\ell$
to that of type $A_\ell$
sending
$$
\F^{\mathcal{D}}( \LD(j)) \simeq \left\{
                                                        \begin{array}{ll}
                                       L(1) \circ L(1)  & \hbox{ if } j=1, \\
                                                          L(j) & \hbox{ otherwise.}
                                                        \end{array}
                                                      \right.
$$
\enEx

We make the following examples for type $B_\ell$ by constructing affinizations directly.

\Ex \label{Ex: B1}
Let $I = \{1,2, \ldots ,\ell \}$ and $\cmA$ the Cartan matrix of type $B_\ell$:
$$ \cmA =
                        \left(
                         \begin{array}{cccccc}
                           2 & -2 & 0  &\cdots & 0 &0\\
                           -1 & 2 & -1  & \cdots  &0&0  \\
                           0 & -1 & 2  & \cdots &  0&0 \\
                           \vdots &  \vdots & \vdots & \ddots & \vdots& \vdots\\
                           0 & 0 &  0 & \cdots & 2 & -1 \\
                           0 & 0 &  0 & \cdots & -1 &2\\
                         \end{array}
                       \right),$$
and $(\alpha_i, \alpha_j)=2\delta(i=j=1)+4\delta(i=j\not=1)-2\delta(|i-j|=1)$ for $i,j\in I$.

Let $R$ be the quiver Hecke algebra associated with $\cmA$ and the parameter $\qQ_{i,j}(u,v)$ defined as follows: for $i,j \in I$ such that $i < j$,
$$ \qQ_{i,j}(u,v) = \left\{
                           \begin{array}{ll}
                             u^2-v & \text{ if } (i,j) = (1,2), \\
                             u-v & \text{ if $j=i+1$ and $(i,j) \ne (1,2)$}  , \\
                             1 & \text{ otherwise. }
                           \end{array}
                         \right.
$$

Let $J = \{1,2,\ldots, \ell-1 \}$ and
$$
\beta_1 = \alpha_1 + \alpha_2,\qquad \beta_j = \alpha_{j+1} \quad (j \in J\setminus \{ 1 \}).
$$
Note that $(\beta_1, \beta_1) = 2$ and $(\beta_j, \beta_j)  = 4$ for $j \ne 1$.
We put
\begin{align*}
\Delta(\beta_1) = L(1,2)_{\mathsf{z}_1}, \qquad \Delta(\beta_j) = L(j+1)_{\mathsf{z}_j} \quad (j \ne 1),
\end{align*}
where the $R(\beta_1)$-module
$L(1,2)_{\mathsf{z}_1} \seteq\bR[\mathsf{z}_1]v$ is defined by
\begin{align*}
e(\nu)v = \delta_{\nu, (1,2)}v,\quad x_j v = \mathsf{z}_1^{(\al_j,\al_j)/2}v, \quad
\tau_1 v = 0.
\end{align*}
Note that $\Delta(\beta_j)$'s are root modules and
$\deg(\mathsf{z}_j)$ is $2$ or $4$ according that
$j=1$ or not.
For $j,k\in J$ with $j\ne k$, we define
$$
\mathsf{R}_{j,k} := R_{ \Delta(\beta_j),\Delta(\beta_k)}
\in \Hom_{R(\beta_j + \beta_k)}( \Delta(\beta_j)\circ \Delta(\beta_k), \Delta(\beta_k)\circ \Delta(\beta_j)  ),
$$
that is,
$$
\mathsf{R}_{j,k} (p \otimes q) = \left\{
                                                      \begin{array}{ll}
                                                        \tau_2 \tau_1 (q \otimes p) & \hbox{ if } j=1, \\
                                                        \tau_1 \tau_2 (q \otimes p) & \hbox{ if } k=1, \\
                                                        \tau_1 (q \otimes p) & \hbox{ otherwise}
                                                      \end{array}
                                                    \right.
$$
for $p\otimes q \in \Delta(\beta_j) \otimes_\bR \Delta(\beta_k) $.
For $j,k \in J$ with $j < k$, we have
\begin{align*}
\mathsf{R}_{k,j}\mathsf{R}_{j,k} =
\left\{
  \begin{array}{ll}
    \mathsf{z}_j^2 - \mathsf{z}_k & \hbox{ if } (j,k) = (1,2), \\
    \mathsf{z}_j- \mathsf{z}_k & \hbox{ if $k=j+1$ and $(j,k) \ne (1,2)$},  \\
    1 & \hbox{ otherwise.}
  \end{array}
\right.
\end{align*}
Then we have the duality datum $ \mathcal{D} = \{ \Delta(\beta_j), \mathsf{z}_j, \mathsf{R}_{j,k}  \}_{j,k\in J } $
and
$$\cmA^{\mathcal{D}} =
                        \left(
                         \begin{array}{cccccc}
                           2 & -2 & 0  &\cdots & 0 &0\\
                           -1 & 2 & -1  & \cdots  &0&0  \\
                           0 & -1 & 2  & \cdots &  0&0 \\
                           \vdots &  \vdots & \vdots & \ddots & \vdots& \vdots\\
                           0 & 0 &  0 & \cdots & 2 & -1 \\
                           0 & 0 &  0 & \cdots & -1 &2\\
                         \end{array}
                       \right)
$$
which is of type $B_{\ell-1}$. Therefore, we have the functor $ \F^{\mathcal{D}}  $ from the category of modules over quiver Hecke algebra of type $B_{\ell-1}$ to that of type $B_{\ell}$
sending
$$
\F^{\mathcal{D}}( L^{\mathcal{D}}(j)) \simeq \left\{
                                                      \begin{array}{ll}
                                                        L(1,2)  & \hbox{ if } j =1, \\
                                                        L(j+1) & \hbox{ otherwise,}
                                                      \end{array}
                                                    \right.
$$
where $L(1,2) = L(1,2)_{\mathsf{z}_1}/ \mathsf{z}_1 L(1,2)_{\mathsf{z}_1}$.

It is easy to check that $\{ \beta_1, \ldots, \beta_{\ell-1} \}$ satisfies $\eqref{Eq: condiotions on beta}$ and
$\cmA^{\mathcal{D}}$ is equal to the matrix $\overline{\cmA}$ defined by $\eqref{Eq: bar A}$. Thus, Proposition~\ref{prop:injcr} implies that
the functor $\F^{\mathcal{D}}$ categorifies
the injective homomorphism $U^-(\g^{\mathcal{D}}) \simeq U^-(\overline{\g}) \longrightarrow U^-(\g) $ and
 $\F^{\mathcal{D}}$ sends the simple modules to simple modules.
By Theorem \ref{Thm: duality and affinization}, for a simple $R^{\mathcal{D}}$-module $N$,
$$ \F^{\mathcal{D}} (  N \hconv \LD(j) ) \simeq \left\{
                                                         \begin{array}{ll}
                                                           \F^{\mathcal{D}}(N) \hconv L(1,2)  & \hbox{ if } j=1, \\
                                                            \F^{\mathcal{D}}(N) \hconv L( j+1)  & \hbox{ otherwise.}
                                                         \end{array}
                                                       \right.
 $$

\enEx

\vskip 0.5em

\Ex \label{Ex: B2}
We use the same notations $I$, $\cmA$, and $\qQ_{i,j}(u,v)$ as in Example \ref{Ex: B1}.

Let $J = \{1,2,\ldots, \ell-1 \}$ and
$$
\beta_1 = 2\alpha_1 + \alpha_2,\qquad \beta_j = \alpha_{j+1} \quad (j \in J\setminus \{1 \}).
$$
Note that $(\beta_j, \beta_j) = 4$ for all $j\in J$.
We define an $R(\beta_1)$-module structure on $L(1,1,2)_{\mathsf{z}_1} := \bR[\mathsf{z}_1] \otimes_{\bR} ( \bR u \oplus \bR v)$ by
\begin{align*}
&e(\nu)(a \otimes u) = \delta_{\nu, (1,1,2)} a \otimes u,
\qquad \qquad \quad \  e(\nu)(a \otimes v) = \delta_{\nu, (1,1,2)} a \otimes v,
 \\
&   x_j(a \otimes u) = \left\{
                       \begin{array}{ll}
                        -\mathsf{z}_1 a \otimes v & \text{if $j=1$,} \\
                         \mathsf{z}_1 a \otimes v  & \text{if $j=2$,} \\
                          \mathsf{z}_1 a \otimes u &\text{otherwise,}
                       \end{array}
                     \right.
 \quad
x_j(a \otimes v) = \left\{
                   \begin{array}{ll}
                     -  a \otimes u  &\text{if $j=1$,} \\
                     a \otimes u & \text{if $j=2$,}\\
                      \mathsf{z}_1 a \otimes v & \text{otherwise,}
                   \end{array}
                 \right.
\\
&  \tau_k (a \otimes u) = \left\{
                            \begin{array}{ll}
                              a \otimes v &\text{if $k=1$,} \\
                              0 & \text{if $k\not=1$,}
                            \end{array}
                          \right.
\qquad \quad  \tau_k (a \otimes v) = 0 \quad\text{for any $k$.}
\end{align*}
We put
\begin{align*}
\Delta(\beta_1) = L(1,1,2)_{\mathsf{z}_1}, \qquad \Delta(\beta_j) = L(j+1)_{\mathsf{z}_{j}} \quad (j\ne 1).
\end{align*}
Note that $\Delta(\beta_j)$'s are root modules and $\deg(\mathsf{z}_j) = 4$ for $j\in J$.
For $j,k\in J$ with $j\ne k$ and $p\otimes q \in \Delta(\beta_j) \otimes_\bR \Delta(\beta_k) $, we define
$$\mathsf{R}_{j,k} :=  R_{\Delta(\beta_j),\,\Delta(\beta_k)}
 \in \Hom_{R(\beta_j + \beta_k)}( \Delta(\beta_j)\circ \Delta(\beta_k), \Delta(\beta_k)\circ \Delta(\beta_j)  ).$$
Then,
\begin{align*}
\mathsf{R}_{k,j}\mathsf{R}_{j,k} =
\left\{
  \begin{array}{ll}
    \mathsf{z}_j - \mathsf{z}_k & \hbox{ if $k=j+1$,} \\
    1 & \hbox{ otherwise,}
  \end{array}
\right.
\end{align*}
for $j,k \in J$ with $j < k$.

Then, we have a duality datum $ \mathcal{D} = \{ \Delta(\beta_j), \mathsf{z}_j,  \mathsf{R}_{j,k}  \}_{j,k\in J} $
and
$$\cmA^{\mathcal{D}} =
                        \left(
                         \begin{array}{cccccc}
                           2 & -1 & 0  &\cdots & 0 &0\\
                           -1 & 2 & -1  & \cdots  &0&0  \\
                           0 & -1 & 2  & \cdots &  0&0 \\
                           \vdots &  \vdots & \vdots & \ddots & \vdots& \vdots\\
                           0 & 0 &  0 & \cdots & 2 & -1 \\
                           0 & 0 &  0 & \cdots & -1 &2\\
                         \end{array}
                       \right)
$$
is of type $A_{\ell-1}$. Therefore, we have the quiver Hecke algebra
$R^{\mathcal{D}}$ of type $A_{\ell-1}$ and
the functor $ \F^{\mathcal{D}}  $ between quiver Hecke algebras of type $A_{\ell-1}$ and $B_\ell$.
Moreover,
$$
\F^{\mathcal{D}}(\LD(j)) \simeq \left\{
                                                      \begin{array}{ll}
                                                        L(1,1,2)  & \hbox{ if } j =1, \\
                                                        L(j+1) & \hbox{ otherwise,}
                                                      \end{array}
                                                    \right.
$$
where $L(1,1,2) = L(1,1,2)_{\mathsf{z}_1}/ \mathsf{z}_1 L(1,1,2)_{\mathsf{z}_1}$.

One can easily show that $\{ \beta_1, \ldots, \beta_{\ell-1} \}$ satisfies $\eqref{Eq: condiotions on beta}$ and
$\cmA^{\mathcal{D}}$ is equal to the matrix $\overline{\cmA}$ defined by $\eqref{Eq: bar A}$. Thus,  Proposition~\ref{prop:injcr} implies
that the functor $\F^{\mathcal{D}}$ categorifies
the injective homomorphism $U^-(\overline{\g}) \longrightarrow U^-(\g) $
and $\F^{\mathcal{D}}$ preserves simple modules.
We have
$$ \F^{\mathcal{D}} (  N \hconv \LD(j) ) \simeq \left\{
                                                         \begin{array}{ll}
                                                           \F^{\mathcal{D}}(N) \hconv L(1,1,2)  & \hbox{ if } j=1, \\
                                                            \F^{\mathcal{D}}(N) \hconv L(j+1)  & \hbox{ otherwise,}
                                                         \end{array}
                                                       \right.
 $$
for a simple $R^{\mathcal{D}}$-module $N$ by Theorem \ref{Thm: duality and affinization}.

\enEx

\end{document}